\documentclass{article}
\usepackage{amsfonts}
\usepackage{mathrsfs}

\usepackage{amsmath}
\usepackage{latexsym,amsfonts,amssymb}

\newcommand{\bfv}{{\bf{v}}}
\newcommand{\bfn}{{\bf{n}}}

\newcommand{\bfnu}{\mbox{\boldmath$\nu$}}

\begin{document}

\title{On the Banach manifold of simple domains in the Euclidean space and applications to free
  boundary problems\thanks{This work is supported by the National Natural  Science Foundation of China
  under grant numbers 11571381}}
\author{Shangbin Cui\\[0.2cm]
  {\small School of Mathematics, Sun Yat-Sen University, Guangzhou, Guangdong 510275,}\\
  {\small People's Republic of China. E-mail:\,cuishb@mail.sysu.edu.cn}}
\date{}
 \maketitle

\begin{abstract}
  In this paper we study the Banach manifold made up of simple $C^{m+\mu}$-domains in the Euclidean space $\mathbb{R}$. This
  manifold is merely a topological or a $C^0$ Banach manifold. It does not possess a differentiable structure. We introduce
  the concept of differentiable point in this manifold and prove that it is still possible to introduce the concept of
  tangent vector and tangent space at a differentiable point. Consequent, it is possible to consider differential equations
  in this Banach space. We show how to reduce some important free boundary problems into differential equations in such a
  manifold and then use the abstract result that we established earlier to study these free boundary problems.
\medskip

\medskip

\end{abstract}

\section{Introduction}
\setcounter{equation}{0}

\hskip 2em
  During the past fifty years, free boundary problems have been one of the most active topics of research in the field of
  partial differential equations. Such problems can be roughly divided into two types: evolutionary type and stationary
  type. To motivate the present work let us put stationary type free boundary problems aside and only consider evolutionary
  type free boundary problems. An evolutionary type free boundary problem contains a domain $\Omega(t)$ depending on time
  $t$ which is not a priori known but has to be determined together with the other unknown function or functions
  such as $u(x,t)$ etc. defined for $t\geqslant 0$ and $x\in\overline{\Omega(t)}$ which satisfies one or a system of partial
  differential equations for $x\in\Omega(t)$ and $t>0$. To determine the domain $\Omega(t)$, in addition to ordinary boundary
  value conditions for the unknown function $u(x,t)$ there must be some extra boundary condition to govern the motion of the free
  boundary $\partial\Omega(t)$. Hence, in an evolutionary type free boundary problem in addition to usual unknown variables
  such as $u$ an extra unknown variable $\Omega$ representing the unknown domain in which partial differential equations are
  expressed has to be involved. Since this extra unknown variable $\Omega$ does not vary in either the real field $\mathbb{R}$
  or the complex field $\mathbb{Z}$, the equation governing the motion of $\Omega(t)$ cannot be directly expressed as a differential
  equation for this extra unknown variable. In current literatures this difficulty is overcome by considering the function
  $\rho(x,t)$ defining the free boundary $\partial\Omega(t)$ with respect to a fixed reference hypersurface $S$. That is, we
  first choose a fixed smooth hypersurface $S$ and express $\partial\Omega(t)$ in accordance with this hypersurface $S$ in
  the following form:
$$
  \partial\Omega(t)=\{x+\rho(x,t){\mathbf{n}}(x):x\in S\},
$$
  where ${\mathbf{n}}$ denotes the normal field of $S$. In this way the non-real unknown variable $\Omega$ is represented
  by a real unknown variable $\rho$, and the equation governing the movement of the free boundary $\partial\Omega(t)$ is
  usually a partial differential equation on the unknown function $\rho(x,t)$.

  The above approach of dealing with evolutionary type free boundary problems has attained great success. But this method
  has the following evident drawback: It can only treat the free boundary problem locally and generally cannot be used to
  study free boundary problems globally; namely, only when $\partial\Omega(t)$ varies in a small neighborhood of the
  reference hypersurface $S$ this method is successful, and when $\partial\Omega(t)$ moves to a position far from $S$ then
  it will no longer be effective. To overcome this difficulty, nowadays an alternative approach has entered the eyes of
  researchers: Take $\Omega$ itself as an unknown variable and let it vary in the set of all bounded domains (with certain
  smoothness) in $\mathbb{R}^n$. For instance, in 1982 when Hamilton \cite{Ham1} studied a free boundary problem
  he considered the unknown domain $\Omega$ as a variable varying in the Frech\`{e}t manifold of all bounded smooth or
  $C^{\infty}$-domains in $\mathbb{R}^n$.

  However, just as Frech\`{e}t space is not as convenient to use as Banach space, Frech\`{e}t manifold is not as convenient
  to use as Banach  manifold. Hence, instead of the manifold $\mathbb{M}_{\infty}$ of smooth bounded domains in $\mathbb{R}^n$,
  one may expect to use the manifold of bounded domains of finite order smoothness as the region for the variable $\Omega$ to
  vary. Unfortunately, such manifold is merely a topological Banach manifold and does not possess a differentiable structure.
  Nondifferentiability of such manifold makes it impossible to perform differential calculus in it in ordinary sense. This is
  perhaps the reason that this approach has not been in the common scope of free boundary problem researchers. Recently, some
  researchers such as the authors of \cite{PruSim} have recognized that this difficulty can be overcome by considering special
  points in such a manifold. To explain this point a little more clearly let us consider the following example: For a positive
  integer $k$ let $\mathbb{M}_k$ be the Banach manifold of all $C^k$-class closed hypersurfaces in $\mathbb{R}^n$. Then every
  $\Sigma\in\mathbb{M}_{k+1}$ can be regarded as a {\em differentiable point} of $\mathbb{M}_k$, and for every such $\Sigma$
  the tangent space $T_{\Sigma}(\mathbb{M}_k)$ can be well-defined because the normal field ${\mathbf{n}}$ of $\Sigma$ is of
  $C^k$-class; see Section 5 of \cite{PruSim}. However, since it seems that the authors of \cite{PruSim} were mainly interested
  in geometric properties of each point $\Sigma$ in $\mathbb{M}_k$, they did not make a systematic study to this manifold.

  The purpose of this paper is to make an elementary study to the Banach manifold $\mathfrak{D}^{m+\mu}(\mathbb{R}^n)$
  of all simple $C^{m+\mu}$-domains in $\mathbb{R}^n$, where $m$ is a positive integer and $0\leqslant\mu\leqslant1$. We
  shall give the definitions of {\em differentiable point} of this manifold and {\em tangent vector} and {\em tangent space}
  at such point in standard language of differential topology. By introducing the concept of {\em standard local chart},
  we shall show that tangent space defined here is equivalent to that as defined in \cite{PruSim}. In this way differential
  calculus can be performed in the nondifferentiable Banach manifold $\mathfrak{D}^{m+\mu}(\mathbb{R}^n)$ and, in particular,
  the concept of differential equations in $\mathfrak{D}^{m+\mu}(\mathbb{R}^n)$ makes sense, which is exactly the motivation
  for us to study differential calculus in the nondifferentiable Banach manifold $\mathfrak{D}^{m+\mu}(\mathbb{R}^n)$.
  Furthermore, we shall also study actions of the following three Lie groups to this manifold: Translation group $G_{tl}$,
  dilation group $G_{dl}$, and rotation group $O(n)$. We shall rigourously prove that actions of these Lie groups to
  $\mathfrak{D}^{m+\mu}(\mathbb{R}^n)$ are differentiable at its differentiable points. These calculuses put our work on
  evolutionary type free boundary problems such as that made in \cite{Cui4} on a solid foundation. Finally, we shall show
  through concrete examples how to reduce evolutionary type free boundary problems into differential equations on the
  Banach manifold $\mathfrak{D}^{m+\mu}(\mathbb{R}^n)$ or certain vector bundles on it.

  Organization of the rest part is as follows. In Section 2 we prove two preliminary lemmas which is the basis of all the
  discussion made in this paper. In Section 3 we give the definition of the Banach manifold structure of
  $\mathfrak{D}^{m+\mu}(\mathbb{R}^n)$ in the language of differential topology, and show how to define $C^k$-points
  ($k\geqslant 1$) in this manifold and the tangent space at such points. We shall show that the definition given here
  is actually equivalent to that as mentioned above. In Section 4 we introduce three kinds of Lie group actions to the
  manifold $\mathfrak{D}^{m+\mu}(\mathbb{R}^n)$. In the last section we present several examples how to reduce evolutionary
  type free boundary problems into differential equation in this manifold or its vector bundles.

\section{Some preliminary results}
\setcounter{equation}{0}

\hskip 2em
  In this section, the notations $S$, $S_1$, $S_2$ denote closed smooth hypersurfaces in $\mathbb{R}^n$. We shall regard all
  closed smooth hypersurfaces in $\mathbb{R}^n$ as $n\!-\!1$-dimensional compact smooth manifolds embedded in $\mathbb{R}^n$,
  and for a such hypersurface $S$, we denote by $\mbox{Chart}(S)$ the set of all local charts of $S$, i.e.,
$$
  \mbox{Chart}(S)=\bigcup_{x\in S}\mbox{Chart}_x(S),
$$
  where $\mbox{Chart}_x(S)$ denotes the set of all local charts of $S$ at the point $x$. Recall that a local chart of $S$ at a
  point $x\in S$ is a couple $(O,\varphi)$ with $O$ an open neighborhood of $x$ in $S$ and $\varphi$ a bijective
  $C^{\infty}$-mapping from $O$ onto the open set $\varphi(O)\subseteq\mathbb{R}^{n-1}$ such that for any $y\in O$,
  the tangent mapping $D_y\varphi:T_y(S)\to\mathbb{R}^{n-1}$ is a linear isomorphism, i.e., $\mbox{rank}D_y\varphi=n\!-\!1$.
  Given a closed smooth hypersurface $S$ in $\mathbb{R}^n$, a nonnegative integer $m$ and a number $\mu\in [0,1]$, the
  notation $C^{m+\mu}(S)$ or alternatively $C^{m+\mu}(S,\mathbb{R})$ denotes the Banach space of all real-valued
  $C^{m+\mu}$-functions on $S$. Recall that a function $f:S\to\mathbb{R}$ is said to be a $C^{m+\mu}$-function on $S$ if for
  any local chart $(O,\varphi)\in\mbox{Chart}(S)$, the function $f\circ\varphi^{-1}:\varphi(O)\to\mathbb{R}$ belongs to
  $C^{m+\mu}(\overline{\varphi(O)})$. Note that for $\mu=1$, the notations $C^{m+\mu}$, $C^{m+\mu}(S)$,
  $C^{m+\mu}(\overline{\varphi(O)})$ indicate $C^{m+1-0}$, $C^{m+1-0}(S)$, $C^{m+1-0}(\overline{\varphi(O)})$, respectively.
  Similarly, given two closed smooth hypersurfaces $S_1$, $S_2$ in $\mathbb{R}^n$, a nonnegative integer $m$ and a number
  $\mu\in [0,1]$, the notation $C^{m+\mu}(S_1,S_2)$ denotes the Banach space of all $C^{m+\mu}$-mappings from $S_1$ to $S_2$.
  Recall that a mapping $F:S_1\to S_2$ is said to be a $C^{m+\mu}$-mapping if for any local charts $(O,\varphi)\in\mbox{Chart}(S_1)$
  and $(U,\psi)\in\mbox{Chart}(S_2)$, the $n\!-\!1$-vector function $\psi\circ F\circ\varphi^{-1}:\varphi(O)\to\psi(U)\subseteq
  \mathbb{R}^{n-1}$ belongs to $C^{m+\mu}(\overline{\varphi(O)},\mathbb{R}^{n-1})$. Note that for $\mu=1$, the notation
  $C^{m+\mu}(S_1,S_2)$ refers to $C^{m+1-0}(S_1,S_2)$. It is clear that for a mapping $F:S_1\to S_2$, $F\in C^{m+\mu}(S_1,S_2)$
  if and only if as a $n$-vector function, $F\in C^{m+\mu}(S_1,\mathbb{R}^n)$.

  Later on, we shall use the notation $\dot{C}^{m+\mu}$ to denote little H\"{o}lder class. For instance, $\dot{C}^{m+\mu}(S)$
  or alternatively $\dot{C}^{m+\mu}(S,\mathbb{R})$ denotes the Banach space of all real-valued ${m+\mu}$-th order little H\"{o}lder
  class functions on $S$, which is defined to be the closure of $C^{\infty}(S)$ in $C^{m+\mu}(S)$. Similarly,
  $\dot{C}^{m+\mu}(S_1,S_2)$ denotes the closure of $C^{\infty}(S_1,S_2)$ in $C^{m+\mu}(S_1,S_2)$.

  Let $X$ and $Y$ be two Banach spaces. Recall that for a positive integer $k$, the notation $L^k(X,Y)$ denotes the Banach
  space of all bounded $k$-linear mappings from $X\times X\times\cdots\times X$ ($k$ times) to $Y$. Recall that a $k$-linear
  mapping $A$ from $X\times X\times\cdots\times X$ ($k$ times) to $Y$ is said to be {\em bounded} if there exists a positive
  constant $C$ such that the following relation holds for all $(x_1,x_2,\cdots,x_k)\in\underbrace{X\times X\times\cdots\times X}_{k}$:
$$
  \|A(x_1,x_2,\cdots,x_k)\|_Y\leqslant C\|x_1\|_X\|x_2\|_X\cdots,\|x_k\|_X.
$$
  The infimum of all such constant $C$ is called the norm of $A$ and is denoted as $\|A\|_{L^k(X,Y)}$. Note that
$$
  L^k(X,Y)\approxeq L(X,L(X,\cdots,L(X,Y)\cdots))\;\; \mbox{($k$-times $L$ and $X$)}.
$$
  Recall that an element $A\in L^k(X,Y)$ is said to be {\em symmetric} if for any permutation $i_1,i_2,\cdots,i_k$ of
  $1,2,\cdots,k$ there holds
$$
  A(x_1,x_2,\cdots,x_k)=A(x_{i_1},x_{i_2},\cdots,x_{i_k}), \quad \forall (x_1,x_2,\cdots,x_k)\in
  \underbrace{X\times X\times\cdots\times X}_{k}.
$$
  The notation $L_s^k(X,Y)$ denotes the Banach subspace of $L^k(X,Y)$ consisting of all bounded symmetric $k$-linear mappings
  from $X\times X\times\cdots\times X$ ($k$ times) to $Y$.

  Let $X$, $Y$ be two Banach spaces. For an open subset $U$ of $X$, a subset $E$ of $Y$ and a nonnegative integer $k$, we use
  the notation $C^k(U,E)$ to denote the set of mappings $F:U\to E$ which are $k$-th order differentiable at every point
  $x\in U$. Here for $1\leqslant j\leqslant k$ we briefly review the $j$-th order derivative $F^{(j)}$ or $D^jF$  of $F$
  as follows. Recall that for $j=1$, $F\in C^1(U,Y)$ if and only if there exists $G_1\in C(U,L(X,Y))$ such that
$$
  \lim_{\|x-x_0\|_X\to 0}\frac{\|F(x)-F(x_0)-G_1(x_0)(x-x_0)\|_Y}{\|x-x_0\|_X}=0,
  \quad \forall x_0\in U,
$$
  and in this case $F'(x)=G_1(x)$, $\forall x\in U$. Inductively, for any integer $2\leqslant j\leqslant k$, if we have already
  known that $F\in C^{j-1}(U,Y)$ and the $j\!-\!1$-th order derivative $F^{(j-1)}$ or $D^{j-1}F$ of $F$ has been defined
  (so that $F^{(j-1)}\in C(U,L^{j-1}(X,Y))$), then $F\in C^j(U,Y)$ if and only if there exists $G_j\in C(U,L^j(X,Y))=
  C(U,L(X,L^{j-1}(X,Y)))$ such that
$$
  \lim_{\|x-x_0\|_X\to 0}\frac{\|F^{(j-1)}(x)-F^{(j-1)}(x_0)-G_j(x_0)(x-x_0)\|_{L^{j-1}(X,Y)}}{\|x-x_0\|_X}=0,
  \quad \forall x_0\in U,
$$
  and in this case $F^{(j)}(x)=G_j(x)$, $\forall x\in U$. It is well-known that $F^{(j)}(x)\in L_s^j(X,Y)$, $\forall x\in U$.

  Let $X$, $Y$ be two Banach spaces and $U$ a bounded open subset of $X$. Given a number $0<\mu\leqslant 1$, we denote by
  $C^{\mu}(\overline{U},Y)$ ($C^{1-0}(\overline{U},Y)$ for $\mu=1$) the set of all mappings $F:U\to Y$ satisfying the following
  conditions: $F$ is bounded on $U$, and there exists a positive constant $C$ such that
$$
  \|F(x_1)-F(x_2)\|_Y\leqslant C\|x_1-x_2\|_X^{\mu}, \quad \forall x_1,x_2\in U.
$$
  The smallest constant $C$ is called the $\mu$-th H\"{o}lder module of $F$ on $U$ and is denoted as $[F]_{\mu;U}$. It is
  clear that $C^{\mu}(\overline{U},Y)$ is a Banach space with norm
$$
  \|F\|_{C^{\mu}(\overline{U},Y)}=\sup_{x\in U}\|F(x)\|_Y+[F]_{\mu;U}, \quad \forall F\in C^{\mu}(\overline{U},Y).
$$
  For a nonnegative integer $k$, we denote by $C^k(\overline{U},Y)$ the set of all $F\in C^k(U,E)$ such that for any $0\leqslant
  j\leqslant k$, $F^{(j)}$ has a bounded continuous extension to $\overline{U}$. It is clear that $C^{k}(\overline{U},Y)$
  is a Banach space with norm
$$
  \|F\|_{C^{k}(\overline{U},Y)}=\sum_{j=0}^k\sup_{x\in U}\|F^{(j)}(x)\|_Y., \quad \forall F\in C^{k}(\overline{U},Y).
$$
  For a positive integer $k$ and a number $0<\mu\leqslant 1$,  we denote by $C^{k+\mu}(\overline{U},Y)$ the set of all
  $F\in C^k(\overline{U},Y)$ such that for any $0\leqslant j\leqslant k$, $F^{(j)}\in C^{\mu}(\overline{U},Y)$. It is clear
  that $C^{k+\mu}(\overline{U},Y)$ is a Banach space with norm
$$
  \|F\|_{C^{k+\mu}(\overline{U},Y)}=\|F\|_{C^{k}(\overline{U},Y)}+\sum_{j=0}^k\|F^{(j)}\|_{C^{\mu}(\overline{U},Y)},
  \quad \forall F\in C^{k+\mu}(\overline{U},Y).
$$
  For $\mu=1$, $C^{k+\mu}(\overline{U},Y)$, $C^{\mu}(\overline{U},Y)$ refer to $C^{k+1-0}(\overline{U},Y)$,
  $C^{1-0}(\overline{U},Y)$, respectively, and for the corresponding norms. In the case $Y=\mathbb{R}$, we often simply
  write $C^{\mu}(\overline{U},\mathbb{R})$, $C^k(\overline{U},\mathbb{R})$, $C^{k+\mu}(\overline{U},\mathbb{R})$ as
  $C^{\mu}(\overline{U})$, $C^k(\overline{U})$, $C^{k+\mu}(\overline{U})$, respectively.

  Finally, let $X$, $Y$ be two Banach spaces and $S$, $S_1$, $S_2$ three closed smooth hypersurfaces in $\mathbb{R}^n$.
  For a nonnegative integer $k$ and a number $0\leqslant\mu\leqslant 1$, the notation $C^{k+\mu}(X,Y)$ denotes the set of
  all mappings $F:X\to Y$ such that for any bounded open set $U\subseteq X$, the restriction of $F$ on $U$ belongs to
  $C^{k+\mu}(\overline{U},Y)$. The notations $C^{k+\mu}(X,S)$,  $C^{k+\mu}(S,Y)$, $C^{k+\mu}(X\times S,Y)$,
  $C^{k+\mu}(X\times S_1,S_2)$ and etc. are defined similarly. We omit the details.
\medskip

  {\bf Lemma 2.1}\ \ {\em Let $S$ be a closed smooth hypersurface in $\mathbb{R}^n$. Let $k$ be a nonnegative integer
  and $0\leqslant\mu\leqslant 1$. Assume that either $k\geqslant 1$ or $k=0$ and $\mu=1$. Then the following relation
  holds:}
\begin{equation}
  [(\rho,x)\mapsto\rho(x)]\in C^{k+\mu}(C^{k+\mu}(S)\times S):=C^{k+\mu}(C^{k+\mu}(S)\times S,\mathbb{R}).
\end{equation}

  {\em Proof}:\ \ By definition, $[(\rho,x)\mapsto\rho(x)]\in C^{k+\mu}(C^{k+\mu}(S)\times S)$ if and only if for any
  local chart $(O,\varphi)\in\mbox{Chart}S$, there holds
$$
  [(\rho,u)\mapsto(\rho\circ\varphi^{-1})(u)]\in C^{k+\mu}(C^{k+\mu}(S)\times\overline{\varphi(O)})
  :=C^{k+\mu}(C^{k+\mu}(S_1)\times\overline{\varphi(O)},\mathbb{R}).
$$
  This will follow if we prove that for any bounded open set $O\subseteq\mathbb{R}^n$, the following relation holds:
\begin{equation}
  [(\eta,u)\mapsto\eta(u)]\in C^{k+\mu}(C^{k+\mu}(\overline{O})\times\overline{O})
  :=C^{k+\mu}(C^{k+\mu}(\overline{O})\times\overline{O},\mathbb{R}).
\end{equation}
  In what follows we denote $f(\eta,u)=\eta(u)$ for $(\eta,u)\in C^{k+\mu}(\overline{O})\times O$, so that $f:
  C^{k+\mu}(\overline{O})\times O\to\mathbb{R}$. We want to prove $f\in C^{k+\mu}(C^{k+\mu}(\overline{O})\times\overline{O})$.
  We fulfill the proof step by step.

  {\em Step 1}:\ \ Consider first the case $k=0$, $\mu=1$. For any $\eta_1,\eta_2\in C^{1-0}(\overline{O})$ and $u_1,u_2\in O$,
  we have
\begin{align*}
  |f(\eta_1,u_1)-f(\eta_2,u_2)|=&|\eta_1(u_1)-\eta_2(u_2)|
\\
  \leqslant&|\eta_1(u_1)-\eta_1(u_2)|+|\eta_1(u_2)-\eta_2(u_2)|
\\
  \leqslant&[\eta_1]_{1,O}\|u_1-u_2\|_{\mathbb{R}^n}+\sup_{u\in O}|\eta_1(u)-\eta_2(u)|.
\end{align*}
  Hence for any bounded open set $\mathcal{U}$ of $C^{1-0}(\overline{O})$ we have
$$
  \sup_{\eta_1,\eta_2\in\mathcal{U}\atop u_1,u_2\in O}|f(\eta_1,u_1)-f(\eta_2,u_2)|
  \leqslant [1+\sup_{\eta\in\mathcal{U}}\|\eta\|_{C^{1-0}(\overline{O})}][\|u_1-u_2\|_{\mathbb{R}^n}+\|\eta_1-\eta_2\|_{C^{1-0}(\overline{O})}].
$$
  This implies that
$$
  [f]_{1;\mathcal{U}\times O}\leqslant 1+\sup_{\eta\in\mathcal{U}}\|\eta\|_{C^{1-0}(\overline{O})}.
$$
  Moreover, we also have
\begin{align*}
  \sup_{\eta\in\mathcal{U},u\in O}|f(\eta,u)|=&\sup_{\eta\in\mathcal{U},u\in O}|\eta(u)|
  \leqslant\sup_{\eta\in\mathcal{U}}\|\eta\|_{C^{1-0}(\overline{O})}.
\end{align*}
  Hence $f\in C^{1-0}(C^{1-0}(\overline{O})\times\overline{O})$.

  {\em Step 2}:\ \ Consider next the case $k=1$, $0\leqslant\mu\leqslant1$. In this case, it is easy to check that the following
  relation holds:
\begin{equation}
  D_{u}f(\eta,u)=\eta'(u), \quad \forall\eta\in C^{1+\mu}(\overline{O}), \;\; \forall u\in O.
\end{equation}
  Moreover, denoting by $l_u$ (for given $u\in O$) the bounded linear functional on $C^{1+\mu}(\overline{O})$ (i.e., $l_u\in
  [C^{1+\mu}(\overline{O})]^*$) defined by
$$
  l_u(h)=h(u), \quad \forall h\in C^{1+\mu}(\overline{O}),
$$
  and by $L$ the continuous mapping from $O$ to $[C^{1+\mu}(\overline{O})]^*=L(C^{1+\mu}(\overline{O}),\mathbb{R})$ defined by
$$
  L(u)=l_u, \quad \forall u\in O,
$$
  we also have
\begin{equation}
  D_{\eta}f(\eta,u)h=L(u)h, \quad \forall\eta\in C^{1+\mu}(\overline{O}), \;\; \forall u\in O.
\end{equation}
  Clearly, if $0<\mu\leqslant 1$ then
\begin{align*}
  \|D_{u}f(\eta_1,u_1)-D_{u}f(\eta_2,u_2)\|_{L(\mathbb{R}^n,\mathbb{R})}=&\;\|\eta_1'(u_1)-\eta_2'(u_2)\|_{L(\mathbb{R}^n,\mathbb{R})}
\\
  \leqslant&\;\|\eta_1'(u_1)-\eta_1'(u_2)\|_{L(\mathbb{R}^n,\mathbb{R})}+\|\eta_1'(u_2)-\eta_2'(u_2)\|_{L(\mathbb{R}^n,\mathbb{R})}
\\
  \leqslant&\;[\eta_1']_{\mu,O}\|u_1-u_2\|_{\mathbb{R}^n}^{\mu}+\sup_{u\in O}\|\eta_1'(u)-\eta_2'(u)\|_{L(\mathbb{R}^n,\mathbb{R})}
\\
  \leqslant&\;[\eta_1']_{\mu,O}\|u_1-u_2\|_{\mathbb{R}^n}^{\mu}+\|\eta_1-\eta_2\|_{C^{1+\mu}(\overline{O})},
\end{align*}
  and
\begin{align*}
  \|D_{u}f(\eta,u)\|_{L(\mathbb{R}^n,\mathbb{R})}=&\;\|\eta'(u)\|_{L(\mathbb{R}^n,\mathbb{R})}
  \leqslant\sup_{u\in O}\|\eta'(u)\|_{L(\mathbb{R}^n,\mathbb{R})}
  \leqslant\|\eta\|_{C^{1+\mu}(\overline{O})}.
\end{align*}
  Hence $D_{u}f\in C^{\mu}(C^{1+\mu}(\overline{O})\times \overline{O},L(\mathbb{R}^n,\mathbb{R}))$ if $0<\mu\leqslant 1$. For $\mu=0$,
  besides the above inequality we have
\begin{align*}
  \|D_{u}f(\eta_1,u_1)-D_{u}f(\eta_2,u_2)\|_{L(\mathbb{R}^n,\mathbb{R})}=&\;\|\eta_1'(u_1)-\eta_2'(u_2)\|_{L(\mathbb{R}^n,\mathbb{R})}
\\
  \leqslant&\;\|\eta_1'(u_1)-\eta_1'(u_2)\|_{L(\mathbb{R}^n,\mathbb{R})}+\|\eta_1-\eta_2\|_{C^{1}(\overline{O})},
\end{align*}
  which implies $\displaystyle\lim_{u_2\to u_1\atop\eta_2\to\eta_1}D_{u}f(\eta_2,u_2)=D_{u}f(\eta_1,u_1)$. Hence
  $D_{u}f\in C(C^{1}(\overline{O})\times \overline{O},L(\mathbb{R}^n,\mathbb{R}))$. Next, it is easy to show
  that $L:O\to[C^{1+\mu}(\overline{O})]^*$ is continuous and, in fact, $L\in C^{1-0}(\overline{O},[C^{1+\mu}(\overline{O})]^*)$.
  Indeed, first we have
\begin{align*}
  \|L(u_1)-L(u_2)\|_{[C^{1+\mu}(\overline{O})]^*}=&\;\sup_{h\in C^{1+\mu}(\overline{O})\atop h\neq 0}
  \frac{|h(u_1)-h(u_2)|}{\|h\|_{C^{1+\mu}(\overline{O})}}
\\
  \leqslant&\;\sup_{h\in C^{1+\mu}(\overline{O})\atop h\neq 0}
  \frac{\sup_{u\in O}\|h'(u)\|_{L(\mathbb{R}^n,\mathbb{R})}}{\|h\|_{C^{1+\mu}(\overline{O})}}\cdot\|u_1-u_2\|_{\mathbb{R}^n}
\\
  \leqslant&\;\|u_1-u_2\|_{\mathbb{R}^n}, \quad \forall u_1,u_2\in O.
\end{align*}
  This shows that $L:O\to[C^{1+\mu}(\overline{O})]^*$ is Lipschitz continuous, i.e., $L\in C^{1-0}(\overline{O},[C^{1+\mu}(\overline{O})]^*)$
  ($0\leqslant\mu\leqslant 1$). Since $D_{\eta}f(\eta,u)=L(u)$ for all $\eta\in C^{1+\mu}(\overline{O})$ and $u\in O$, it follows that
  $D_{\eta}f\in C^{1-0}(C^{1+\mu}(\overline{O})\times\overline{O},[C^{1+\mu}(\overline{O})]^*)\subseteq
  C^{\mu}(C^{1+\mu}(\overline{O})\times\overline{O},L(C^{1+\mu}(\overline{O}),\mathbb{R}))$ ($0\leqslant\mu\leqslant 1$).
  Therefore, $f\in C^{1+\mu}(C^{1+\mu}(\overline{O})\times\overline{O})$ ($0\leqslant\mu\leqslant 1$).

  {\em Step 3}:\ \ For $k=2$ we have
\begin{align*}
  &D_{u^2}^2f(\eta,u)=\eta''(u), \quad \forall\eta\in C^{2}(\overline{O}), \;\; \forall u\in O,
\\
  &D_{u\eta}^2f(\eta,u)h=h'(u), \quad \forall\eta,h\in C^{2}(\overline{O}), \;\; \forall u\in O,
\\
  &D_{\eta^2}^2f(\eta,u)=0, \quad \forall\eta\in C^{2}(\overline{O}), \;\; \forall u\in O.
\end{align*}
  In general, for $k\geqslant 2$ we have
\begin{align*}
  &D_{u^k}^kf(\eta,u)=\eta^{(k)}(u), \quad \forall\eta\in C^{k}(\overline{O}), \;\; \forall u\in O,
\\
  &D_{u^{k-1}\eta}^kf(\eta,u)h=h^{(k-1)}(u), \quad \forall\eta,h\in C^{k}(\overline{O}), \;\; \forall u\in O,
\\
  &D_{u^{k-j}\eta^j}^kf(\eta,u)=0, \quad \forall\eta\in C^{k}(\overline{O}), \;\; \forall u\in O, \;\; j=2,3,\cdots,k.
\end{align*}
  Using these relations and some similar arguments as in Step 2, we see that the desired assertion is also true for
  the general case $k\geqslant 2$, $0\leqslant\mu\leqslant1$. We omit the details here. $\quad\Box$
\medskip

  {\bf Lemma 2.2}\ \ {\em Let $S_1,S_2$ be two closed smooth hypersurfaces in $\mathbb{R}^n$. Let $k$ be a positive integer
  and $0\leqslant\mu\leqslant 1$. Let $\chi\in C^{k+\mu}(C^{k+\mu}(S_1)\times S_2,S_1)$. Then the following relation
  holds:}
\begin{equation}
  [(\rho,y)\mapsto\rho(\chi(\rho,y))]\in C^{k+\mu}(C^{k+\mu}(S_1)\times S_2):=C^{k+\mu}(C^{k+\mu}(S_1)\times S_2,\mathbb{R}).
\end{equation}

  {\em Proof}:\ \ Let $\kappa:C^{k+\mu}(S_1)\times S_1\to\mathbb{R}$ be the mapping defined by
$$
  \kappa(\rho,x)=\rho(x), \quad \forall\rho\in C^{k+\mu}(S_1),\;\; \forall x\in S_1.
$$
  Then we have
$$
  \rho(\chi(\rho,y))=\kappa(\rho,\chi(\rho,y)), \quad \forall\rho\in C^{k+\mu}(S_1),\;\; \forall y\in S_2.
$$
  Since by Lemma 2.1 we have $\kappa\in C^{k+\mu}(C^{k+\mu}(S_1)\times S_1,\mathbb{R})$, and the given condition guarantees
  that $\chi\in C^{k+\mu}(C^{k+\mu}(S_1)\times S_2,S_1)$, it follows immediately that
$$
  [(\rho,y)\mapsto\rho(\chi(\rho,y))]=[(\rho,y)\mapsto\kappa(\rho,\chi(\rho,y))]\in
  C^{k+\mu}(C^{k+\mu}(S_1)\times S_2,\mathbb{R}).
$$
  Here we have used the following simple proposition: Let $X,Y,Z$ be Banach spaces, and $S_1,S_2$ closed smooth hypersurfaces
  in $\mathbb{R}^n$. Let $k$ be a positive integer and $0\leqslant\mu\leqslant 1$. Assume that $f\in C^{k+\mu}(X\times S_1,Y)$,
  $g\in C^{k+\mu}(X\times S_1,S_2)$ and $h\in C^{k+\mu}(Y\times S_2,Z)$. Then
$$
  [(\rho,x)\mapsto h(f(\rho,y),g(\rho,y))]\in C^{k+\mu}(X\times S_1,Z)
$$
  (note that this is false for $k=0$ and $0<\mu<1$). This proves the desired assertion. $\quad\Box$
\medskip

  {\bf Corollary 2.3}\ \ {\em Let $S$ be a closed smooth hypersurface in $\mathbb{R}^n$. Let $k$ be a positive integer
  and $0\leqslant\mu\leqslant 1$. Let $f\in C^{k+1}(S\times\mathbb{R},S)$. Then the following relation
  holds:}
\begin{equation}
  [(\rho,x)\mapsto\rho(f(x,\rho(x)))]\in C^{k+\mu}(C^{k+\mu}(S)\times S):=C^{k+\mu}(C^{k+\mu}(S)\times S,\mathbb{R}).
\end{equation}

  {\em Proof}:\ \ Let $\chi: C^{k+\mu}(S)\times S\to S$ be the map defined by $\chi(\rho,x)=f(x,\rho(x))$, $\forall
  (\rho,x)\in C^{k+\mu}(S)\times S$. Then it is not hard to verify that $\chi\in C^{k+\mu}(C^{k+\mu}(S)\times S,S)$.
  Indeed, to verify this assertion we only need to prove that for any two local charts $(U,\varphi)$ and $(V,\psi)$ of
  $S$, where $U,V$ are two open sets in $S$ and $\varphi:U\to\varphi(U)\subseteq\mathbb{R}^{n-1}$, $\psi:V\to\psi(V)
  \subseteq\mathbb{R}^{n-1}$ are two $C^{k+\mu}$-diffeomorphisms, the map $(\xi,u)\mapsto\tilde{\chi}(\xi,u)=
  \tilde{f}(u,\xi(u))$, $\forall(\xi,u)\in C^{k+\mu}(\overline{\varphi(U)})\times\varphi(U)$, where $\tilde{f}(u,s)=
  \psi(f(\varphi^{-1}(u),s))$, $\forall(u,s)\in\varphi(U)\times\mathbb{R}$, belongs to $C^{k+\mu}(C^{k+\mu}(\overline{\varphi(U)})
  \times\varphi(U),\mathbb{R}^{n-1})$. Since $f\in C^{k+1}(S\times\mathbb{R},S)$, we have $\tilde{f}\in
  C^{k+1}(\overline{\varphi(U)}\times\mathbb{R},\mathbb{R}^{n-1})$. It easily follows that $\tilde{\chi}\in
  C^{k+\mu}(C^{k+\mu}(\overline{\varphi(U)})\times\varphi(U),\mathbb{R}^{n-1})$. Hence by lemma 2.2, the desired
  assertion follows. $\quad\Box$
\medskip

  Let $X$ and $X_0$ be two Banach spaces such that $X_0$ is a densely embedded Banach subspace of $X$. Let $U_0$ be an open subset
  of $X_0$. Let $Y$ be another Banach space. For a map $F:U_0\to Y$ and a point $x_0\in U_0$, we say $x_0$ is a {\em $d^k$-point}
  of $F$, where $k$ is a positive integer, if for each $1\leqslant j\leqslant k$ there exists a corresponding operator $A_j
  \in L_s^j(X,Y)$ such that the following relation holds:
$$
  \lim_{\|x-x_0\|_{X_0}\to 0}\frac{\Big\|F(x)-F(x_0)
  -\displaystyle\sum_{j=1}^k\frac{1}{j!}A_j(x-x_0,x-x_0,\cdots,x-x_0)\Big\|_{Y}}{\|x-x_0\|_{X_0}^k}=0.
$$
  $A_j$ is called the $j$-th order differential or $j$-th order Fr\'{e}chet derivative of $F$ at $x_0$ and is denoted as
  $D^jF(x_0)=A_j$ or $F^{(j)}(x_0)=A_j$, $j=1,2,\cdots,k$. In particular, for $j=1,2$ the operators $F^{(1)}(x_0)$ and
  $F^{(2)}(x_0)$ are also denoted as $F'(x_0)$ and $F''(x_0)$, respectively. Note that since $L^j(X,Y)\subseteq
  L^j(X_0,Y)$, it follows that $A_j\in L_s^j(X,Y)$ implies $A_j\in L_s^j(X_0,Y)$. From this fact it can be easily seen
  that if $x_0$ is a $d^k$-point of $F$ then it is also a $d^j$-point of $F$ for any $1\leqslant j\leqslant k\!-\!1$.
\medskip

  {\em Remark}.\ Note that if $X_0$ is not dense in $X$, then the operators $A_1$, $A_2$, $\cdots$, $A_k$ might not be uniquely
  determined by $F$ and $x_0$. The reason is that the above relation does not use values of these operators outside $X_0$,
  so that it is possible to change values of them in any subspace of $X$ which is complementary to $\bar{X}_0$ without
  changing the above relation if $X_0$ is not dense in $X$. If, however, $X_0$ is dense in $X$, then clearly $A_1$, $A_2$,
  $\cdots$, $A_k$ are uniquely determined by $F$ and $x_0$.
\medskip

  Let $X,X_0,Y$ and $U_0$ be as above. Given a positive integer $k$, we use the notation $\mathfrak{C}^k(U_0;X,Y)$ to
  denote the set of all mappings $F:U_0\to Y$ satisfing the following two conditions:
\begin{description}
\item[] (1)\ \ All points in $U_0$ are $d^k$-points of $F$;
\item[] (2)\ \ $[x\mapsto F^{(j)}(x)]\in C(U_0,L^j(X,Y))$, $j=1,2,\cdots,k$, where $U_0$ uses the topology of $X_0$.
\end{description}
  It is clear that $\mathfrak{C}^k(U_0;X,Y)\subseteq C^k(U_0,Y)$, where $C^k(U_0,Y)$ denotes the set of all $k$-th order
  continuously differentiable mappings $F:U_0\subseteq X_0\to Y$, and an element $F\in C^k(U_0,Y)$ belongs to
  $\mathfrak{C}^k(U_0;X,Y)$ if and only if for each $1\leqslant j\leqslant k$ and any $x\in U_0$, $D^jF(x)\in L_s^j(X_0,Y)$
   can be extended into an operator belonging to $L_s^j(X,Y)$, and after extension $D^jF\in C(U_0,L^j(X,Y))$, where $U_0$ uses
   the topology of $X_0$. Hence the condition $F\in\mathfrak{C}^k(U_0;X,Y)$ is stronger than the condition $F\in C^k(U_0,Y)$.
\medskip

  {\bf Lemma 2.4}\ \ {\em Let $S$ be a closed smooth hypersurface in $\mathbb{R}^n$. Let $m$ be a positive integer and
  $0\leqslant\mu\leqslant 1$. Let  $X=\dot{C}^{m+\mu}(S)$ and $X_0=\dot{C}^{m+1+\mu}(S)$. Let $f\in C^{m+1}(S\times\mathbb{R},
  S)$ and $F(\rho)=[\rho\mapsto [x\mapsto\rho(f(x,\rho(x)))]]$, $\rho\in X_0$. Then the following relation
  holds:}
\begin{equation}
  F\in \mathfrak{C}^1(X_0;X,X), \quad \mbox{i.e.}, \quad  F'\in C(X_0,L(X)).
\end{equation}

  {\em Proof}:\ \ We denote $F(\rho)=[\rho\mapsto [x\mapsto\rho(f(x,\rho(x)))]]$. From the proofs of Lemmas 2.1 and 2.2 we see
  that $F\in C^1(C^{m+1+\mu}(S),{C}^{m+\mu}(S))$. What we need to prove is that for any $\rho\in C^{m+1+\mu}(S)$, the map
  $F'(\rho)\in L(X_0,X)$ can be extended into a map in $L(X,X)$, and $[\rho\mapsto F'(\rho)]\in C(X_0,L(X,X))$. Indeed, for
  any $\rho\in X_0$ we have
\begin{equation}
  F'(\rho)\eta=\rho'(f(\cdot,\rho(\cdot)))\partial_{\rho}\!f(\cdot,\rho(\cdot))\eta+\eta(f(\cdot,\rho(\cdot))), \quad
  \forall\eta\in X_0.
\end{equation}
  From this expression we easily see that the operator $F'(\rho)\in L(X_0,X)$ can be extended into a bounded linear
  operator in $X$: $F'(\rho)\in L(X,X)$. Moreover, since $X_0$ is dense in $X$, the extension is unique. This proves
  the desired assertion. $\quad\Box$
\medskip

  {\bf Lemma 2.5}\ \ {\em Let $S$ be a closed smooth hypersurface in $\mathbb{R}^n$. Let $m$ be a positive integer and
  $0\leqslant\mu\leqslant 1$. Let $X=\dot{C}^{m+\mu}(S)$, $X_0=\dot{C}^{m+1+\mu}(S)$ and $X_1=\dot{C}^{m+2+\mu}(S)$.
  Let $f\in C^{m+2}(S\times\mathbb{R},S)$ and $F(\rho)=[\rho\mapsto [x\mapsto\rho(f(x,\rho(x)))]]$, $\rho\in X_1$. Then
  the following relation holds:}
\begin{equation}
  F''\in C(X_1,L^2_s(X_0,X)).
\end{equation}

  {\em Proof}:\ \ A simple computation shows that for any $\rho\in X_0=\dot{C}^{m+2+\mu}(S)$,
\begin{eqnarray}
  F''(\rho)(\eta,\zeta)&=&[x\mapsto\rho''(f(x,\rho(x)))(\partial_{\rho}\!f(x,\rho(x))\eta(x),\zeta(x))
  +\rho'(f(x,\rho(x)))\partial_{\rho}^2\!f(x,\rho(x))(\eta(x),\zeta(x))
\nonumber\\
  &&+\eta'(f(x,\rho(x)))\partial_{\rho}\!f(x,\rho(x))\zeta(x)+\zeta'(f(x,\rho(x)))\partial_{\rho}\!f(x,\rho(x))\eta(x)],
  \quad \forall\eta\in X_1.
\end{eqnarray}
  From this expression we easily see that (2.9) holds. $\quad\Box$
\medskip

\section{Banach manifolds of simple domains in ${\mathbb{R}}^n$}
\setcounter{equation}{0}

\hskip 2em
  In this section we make a basic study to Banach manifolds of simple domains in ${\mathbb{R}}^n$. Such a manifold is only a
  topological manifold and does not possess a differentiable structure. However, some of its embedded Banach submanifolds
  possess differentiable structure in its topology, so that the technique of differential calculus is still possible in
  such Banach submanifolds.

  As usual, given a nonnegative integer $m$, a number $\mu\in [0,1]$, a bounded open set $\Omega\subseteq{\mathbb{R}}^n$ and a
  sufficiently smooth (e.g., smooth up to order $k$ for some integer $k\geqslant m+\mu$) closed hypersurface $S\subseteq{\mathbb{R}}^n$,
  the notations $C^{m+\mu}(\overline{\Omega})$ and $C^{m+\mu}(S)$ denote the usual
  $m\!+\!\mu$-th order H\"{o}lder spaces on $\overline{\Omega}$ and $S$, respectively, and the notation $C^{m+\mu}(\overline{\Omega},{\mathbb{R}}^n)$
  denotes the usual $m\!+\!\mu$-th order $n$-vector H\"{o}lder space on $\overline{\Omega}$, with the cases $\mu=0,1$ understood in conventional sense.
  We use the notation $\dot{C}^{m+\mu}(\overline{\Omega})$ to denote the closure of $C^{\infty}(\overline{\Omega})$ in $C^{m+\mu}(\overline{\Omega})$,
  and similarly for the notations $\dot{C}^{m+\mu}(S)$ and $\dot{C}^{m+\mu}(\overline{\Omega},{\mathbb{R}}^n)$. The last three spaces are called
  $m\!+\!\mu$-th order little H\"{o}lder spaces. A significant difference between little H\"{o}lder spaces and H\"{o}lder spaces is that for
  nonnegative integers $k,m$ and real numbers $\mu,\nu\in [0,1]$, if $k+\nu>m+\mu$ then $\dot{C}^{k+\nu}(\overline{\Omega})$ (resp.
  $\dot{C}^{k+\nu}(S)$, $\dot{C}^{k+\nu}(\overline{\Omega},{\mathbb{R}}^n)$) is dense in $\dot{C}^{m+\mu}(\overline{\Omega})$ (resp.
  $\dot{C}^{m+\mu}(S)$, $\dot{C}^{m+\mu}(\overline{\Omega},{\mathbb{R}}^n)$), but $C^{k+\nu}(\overline{\Omega})$ (resp. $C^{k+\nu}(S)$,
  $C^{k+\nu}(\overline{\Omega},{\mathbb{R}}^n)$) is not dense in $C^{m+\mu}(\overline{\Omega})$ (resp. $C^{m+\mu}(S)$,
  $C^{m+\mu}(\overline{\Omega},{\mathbb{R}}^n)$).
\medskip

  {\bf Definition 3.1}\ \ {\em Let $m$ be a positive integer and $0\leqslant\mu\leqslant 1$. An open set $\Omega\subseteq{\mathbb{R}}^n$ is said to be
  a simple $C^{m+\mu}$-domain if $\Omega$ is $C^{m+\mu}$-diffeomorphic to the open unit sphere ${\mathbb{B}}^n$ in ${\mathbb{R}}^n$, i.e., there exists a bijective
  mapping $\Phi:\overline{{\mathbb{B}}^n}\to\overline{\Omega}$ satisfying the following properties:
$$
  \Phi\in C^{m+\mu}(\overline{{\mathbb{B}}^n},{\mathbb{R}}^n) \quad \mbox{and} \quad \Phi^{-1}\in C^{m+\mu}(\overline{\Omega},{\mathbb{R}}^n).
$$
  We use the notation $\mathfrak{D}^{m+\mu}({\mathbb{R}}^n)$ to denote the set of all simple $C^{m+\mu}$-domains in ${\mathbb{R}}^n$. If instead of
  $C^{m+\mu}$ the notation $\dot{C}^{m+\mu}$ is used in the above relations, then the notation $\dot{\mathfrak{D}}^{m+\mu}({\mathbb{R}}^n)$ is used
  correspondingly.}
\medskip

  From \cite{Ham1} we know that all smooth simple domains in ${\mathbb{R}}^n$ form a Frech\'{e}t manifold $\mathfrak{D}^{\infty}({\mathbb{R}}^n)$
  built on the Frech\'{e}t space $C^{\infty}({\mathbb{S}}^{n-1})$, with tangent space at the point $\Omega\in\mathfrak{D}^{\infty}({\mathbb{R}}^n)$
  being $T_{\Omega}(\mathfrak{D}^{\infty}({\mathbb{R}}^n))=C^{\infty}(\partial\Omega)$. In application, however, just as Frech\'{e}t space
  is not as convenient to use as Banach space, Frech\'{e}t manifold is not as convenient to use as Banach manifold. Hence, in what follows
  we introduce a local chart system for $\mathfrak{D}^{m+\mu}({\mathbb{R}}^n)$ ($m\in \mathbb{N}$, $0\leqslant\mu\leqslant 1$) to make it into
  a Banach manifold.
\medskip

  {\bf Lemma 3.2}\ \ {\em Let $m$ be a positive integer and $0\leqslant\mu\leqslant 1$. Let $\Phi\in
  C^{m+\mu}(\overline{{\mathbb{B}}^n},{\mathbb{R}}^n)$ such that $\Phi:{\mathbb{B}}^n\to\Omega=\Phi({\mathbb{B}}^n)$ is
  a bijection, and $\Phi^{-1}\in C^{m+\mu}(\overline{\Omega},{\mathbb{R}}^n)$. Then there exists $\varepsilon>0$ such that
  for any $\Psi\in C^{m+\mu}(\overline{{\mathbb{B}}^n},{\mathbb{R}}^n)$, if
$$
  \|\Psi-\Phi\|_{C^{1}(\overline{{\mathbb{B}}^n},{\mathbb{R}}^n)}<\varepsilon,
$$
  then $\Psi$ is invertible, and $\Psi^{-1}\in C^{m+\mu}(\overline{Q},{\mathbb{R}}^n)$, where $Q=\Psi({\mathbb{B}}^n)$.}
\medskip

  {\em Proof}.\ \ Let $M=\displaystyle\max_{y\in\overline{\Omega}}\|D\Phi^{-1}(y)\|$. Here $\|\cdot\|$ denotes the norm of
  $n\times n$ matrices. Since $D\Phi(x)^{-1}=D\Phi^{-1}(\Phi(x))$, $x\in\overline{{\mathbb{B}}^n}$, it follows by a standard
  result in matrix theory that for any $x\in\overline{{\mathbb{B}}^n}$, if a $n\times n$ matrix $A$ satisfies
  $\|D\Phi(x)-A\|<1/M$ then $A$ is invertible. Hence, by letting $\varepsilon_1=1/M$, we see that for any $\Psi\in
  C^{m+\mu}(\overline{{\mathbb{B}}^n},{\mathbb{R}}^n)$, if $\|\Psi-\Phi\|_{C^{1}(\overline{{\mathbb{B}}^n},{\mathbb{R}}^n)}
  <\varepsilon_1$ then $\Psi$ is locally invertible, and the local
  inverse $\Psi^{-1}$ belongs to $C^{m+\mu}$-class. Since $\overline{\Omega}=\Phi(\overline{{\mathbb{B}}^n})$ and $\Phi\in
  C^{m+\mu}(\overline{{\mathbb{B}}^n},{\mathbb{R}}^n)$, $m\geqslant1$, there exists constant $C>0$ such that for any two
  points $u,v\in\overline{\Omega}$ there exists corresponding $C^1$-curve in $\overline{\Omega}$ with length not larger than
  $C\|u-v\|$ to connect $u$ and $v$. It follows by the mean value inequality that
$$
  |\Phi^{-1}(u)-\Phi^{-1}(v)|\leqslant CM|u-v|, \quad \forall u,v\in\overline{\Omega},
$$
  which implies
$$
  |x-y|\leqslant CM|\Phi(x)-\Phi(y)|, \quad \forall x,y\in\overline{{\mathbb{B}}^n}.
$$
  From this we can easily deduce that there exists $\varepsilon_2>0$ such that for any $\Psi\in
  C^{m+\mu}(\overline{{\mathbb{B}}^n},{\mathbb{R}}^n)$, if $\|\Psi-\Phi\|_{C^{1}(\overline{{\mathbb{B}}^n},{\mathbb{R}}^n)}
  <\varepsilon_2$ then
$$
  |x-y|\leqslant 2CM|\Psi(x)-\Psi(y)|, \quad \forall x,y\in\overline{{\mathbb{B}}^n},
$$
  which implies that $\Psi$ is invertible. Hence, by letting $\varepsilon=\min\{\varepsilon_1,\varepsilon_2\}$, the
  desired assertion follows. $\quad\Box$
\medskip

  Let $S$ be a closed $C^{m+1+\mu}$-hypersurface in $\mathbb{R}^n$ enclosing a simple $C^{m+1+\mu}$-domain $Q\in
  \mathfrak{D}^{m+1+\mu}({\mathbb{R}}^n)$. Let $\bfn$ be the unit normal field of $S$, outward pointing with respect
  to $Q$. It is known that $\bfn\in C^{m+\mu}(S)$. For $\delta>0$ we denote
\begin{equation}
  \mathcal{R}=\{x\in{\mathbb{R}}^n:d(x,S)<4\delta\}.
\end{equation}
  Clearly, if $\delta>0$ is sufficiently small then the mapping $\Upsilon:S\times(-4\delta,4\delta)\to \mathcal{R}$,
\begin{equation}
   \Upsilon(x,t)=x+t\bfn(x), \quad \forall x\in S,\;\; \forall t\in (-4\delta,4\delta)
\end{equation}
  is a $C^{m+\mu}$-diffeomorphism of $S\times(-4\delta,4\delta)$ onto $\mathcal{R}$. Let $\Pi$ and $\Lambda$ be
  compositions of $\Upsilon^{-1}:\mathcal{R}\to S\times(-4\delta,4\delta)$ with the projections $P_1$ and $P_2$ of
  $S\times (-4\delta,4\delta)$ onto $S$ and $(-4\delta,4\delta)$, respectively, i.e., $\Pi:\mathcal{R}\to S$,
  $\Lambda:\mathcal{R}\to(-4\delta,4\delta)$,
\begin{equation}
  \Pi(y)=P_1(\Upsilon^{-1}(y)), \quad \Lambda(y)=P_2(\Upsilon^{-1}(y)), \quad \forall y\in\mathcal{R}.
\end{equation}
  Clearly,
$$
  \Pi(\Upsilon(x,t))=x, \quad \Lambda(\Upsilon(x,t))=t, \quad \forall x\in S, \;\; \forall t\in(-4\delta,4\delta),
$$
  and
$$
  y=\Pi(y)+\Lambda(y)\bfn(\Pi(y)),\;\; \forall y\in \mathcal{R}.
$$
  It is easy to see that if $\delta$ is small enough which we assume to be true then for any $y\in \mathcal{R}$,
  $\Pi(y)$ is the point in $S$ nearest to $y$ and $\Lambda(y)$ is the algebraic distance of $y$ to $S$, i.e,
$$
  \Lambda(y)=-d(y,S)\;\; \mbox{if}\;\, y\in\mathcal{R}\cap Q \quad \mbox{and} \quad
  \Lambda(y)=d(y,S)\;\; \mbox{if}\;\, y\in\mathcal{R}\backslash Q.
$$
  We denote
\begin{equation}
  \mathcal{O}=\{\rho\in C^{m+\mu}(S):\,\|\rho\|_{C^{1}(S)}<\delta\}.
\end{equation}
  Since $\rho\in\mathcal{O}$ implies $\displaystyle\max_{x\in S}|\rho(x)|<\delta$, it follows that for any
  $\rho\in\mathcal{O}$ the mapping $\theta_{\rho}:S\to{\mathbb{R}}^n$,
\begin{equation}
   \theta_{\rho}(x)=x+\rho(x)\bfn(x), \quad \forall x\in S
\end{equation}
  is a $C^{m+\mu}$-diffeomorphism of $S$ onto its image $S_{\rho}=\theta_{\rho}(S)$, and $S_{\rho}$ is a closed
  $C^{m+\mu}$-hypersurface.  We denote by $\Omega_{\rho}$ the domain enclosed by $S_{\rho}$.
\medskip

  {\bf Lemma 3.3}\ \ {\em Let assumptions and notations be as above. If $\delta>0$ is sufficiently small, then for any
  $\rho\in\mathcal{O}$, $\Omega_{\rho}$ is a simple $C^{m+\mu}$-domain.}
\medskip

  {\em Proof}.\ \ Choose a function $\phi\in C^{\infty}[0,\infty)$ such that
$$
  \phi(t)=0\;\; \mbox{for}\;\, 0\leqslant t\leqslant\frac{1}{2}, \quad
  \phi(t)=1\;\; \mbox{for}\;\, t\geqslant 1, \quad \mbox{and} \quad
  \phi'(t)\geqslant 0\;\; \mbox{for}\;\, t\geqslant 0.
$$
  Let $\Phi\in C^{m+1+\mu}(\overline{{\mathbb{B}}^n},{\mathbb{R}}^n)$ be a $C^{m+1+\mu}$-diffeomorphism of
  $\overline{{\mathbb{B}}^n}$ to $\overline{Q}$, i.e., $\overline{Q}=\Phi(\overline{{\mathbb{B}}^n})$, $\Phi$ is
  a bijection of $\overline{{\mathbb{B}}^n}$ onto $\overline{Q}$, and $\Phi^{-1}\in C^{m+1+\mu}(\overline{Q},{\mathbb{R}}^n)$.
  Given $\rho\in\mathcal{O}$, we define a map $\Psi_{\rho}:\overline{{\mathbb{B}}^n}\to{\mathbb{R}}^n$ as follows:
$$
  \Psi_{\rho}(x)=\Phi(x)+\phi(|x|)\rho(\Phi(\omega))\bfn(\Phi(\omega)), \quad x\in\overline{{\mathbb{B}}^n}.
$$
  where $\omega=x/|x|$ for $x\in\backslash\{0\}$. Clearly, $\Psi_{\rho}\in C^{m+\mu}(\overline{{\mathbb{B}}^n},{\mathbb{R}}^n)$,
  and $\Psi_{\rho}(\overline{{\mathbb{B}}^n})=\overline{\Omega}_{\rho}$. Moreover, it is easy to see that there exists
  constant $C>0$ such that
$$
  \|\Psi_{\rho}-\Phi\|_{C^{1}(\overline{{\mathbb{B}}^n},{\mathbb{R}}^n)}\leqslant C\|\rho\|_{C^{1}(S)}<C\delta.
$$
  It follows by Lemma 3.2 that if $\delta>0$ is sufficiently small, then $\Psi_{\rho}:{\mathbb{B}}^n\to\Omega_{\rho}$ is
  a bijection, and $\Psi_{\rho}^{-1}\in C^{m+\mu}(\overline{\Omega}_{\rho},{\mathbb{R}}^n)$. Hence $\Omega_{\rho}\in
  \mathfrak{D}^{m+\mu}({\mathbb{R}}^n)$. This proves the desired assertion. $\quad\Box$
\medskip

  In the following we assume that $S$ is smooth and denote
\begin{equation}
  \mathcal{U}=\{\Omega_{\rho}:\rho\in\mathcal{O}\}.
\end{equation}
  We define
\begin{equation}
  \varphi:\mathcal{U}\to C^{m+\mu}(S), \quad \varphi(\Omega_{\rho})=\rho, \;\;\forall\rho\in\mathcal{O}.
\end{equation}
  We call the pair $(\mathcal{U},\varphi)$ a {\em regular local chart} of $\mathfrak{D}^{m+\mu}({\mathbb{R}}^n)$, and call the
  closed hypersurface $S$ the {\em base hypersurface} of this local chart (the phrase ``regular'' refers to the fact that the
  base hypersurface $S$ is smooth). Clearly, $C^{m+\mu}(S)\cong C^{m+\mu}({\mathbb{S}}^{n-1})$, i.e., $C^{m+\mu}(S)$ and
  $C^{m+\mu}({\mathbb{S}}^{n-1})$ are isomorphic to each other as Banach spaces. We now denote by $\mathscr{A}$ the set of all
  regular local chart of $\mathfrak{D}^{m+\mu}({\mathbb{R}}^n)$:
\begin{equation}
  \mathscr{A}=\{(\mathcal{U},\varphi):(\mathcal{U},\varphi)\;\mbox{is a regular local chart of $\mathfrak{D}^{m+\mu}({\mathbb{R}}^n)$}\}.
\end{equation}
  Regular local chart of $\dot{\mathfrak{D}}^{m+\mu}({\mathbb{R}}^n)$ is defined similarly, with all $C^{m+\mu}(S)$ above replaced by
  $\dot{C}^{m+\mu}(S)$. We denote by $\dot{\mathscr{A}}$ the set of all regular local chart of $\dot{\mathfrak{D}}^{m+\mu}({\mathbb{R}}^n)$:
\begin{equation}
  \dot{\mathscr{A}}=\{(\dot{\mathcal{U}},\varphi):(\dot{\mathcal{U}},\varphi)\;\mbox{is a regular local chart of
  $\dot{\mathfrak{D}}^{m+\mu}({\mathbb{R}}^n)$}\}.
\end{equation}

  {\em Remark}.\ \ From the definition of regular local chart we see that if $(\mathcal{U},\varphi)$ is a regular
  local chart of $\mathfrak{D}^{m+\mu}({\mathbb{R}}^n)$, then for any $k\in{\mathbf{Z}}_+$, if $\Omega\in\mathcal{U}
  \cap\mathfrak{D}^{m+k+\mu}({\mathbb{R}}^n)$ then $\varphi(\Omega)\in C^{m+k+\mu}(S)$. Similarly, if
  $(\dot{\mathcal{U}},\varphi)$ is a regular local chart of $\dot{\mathfrak{D}}^{m+\mu}({\mathbb{R}}^n)$, then for any
  $k\in{\mathbf{Z}}_+$, if $\Omega\in\dot{\mathcal{U}}\cap\dot{\mathfrak{D}}^{m+k+\mu}({\mathbb{R}}^n)$ then $\varphi(\Omega)\in
  \dot{C}^{m+k+\mu}(S)$.
\medskip

  {\bf Theorem 3.4}\ \ {\em $(\mathfrak{D}^{m+\mu}({\mathbb{R}}^n),\mathscr{A})$ $($resp.
  $(\dot{\mathfrak{D}}^{m+\mu}({\mathbb{R}}^n),\dot{\mathscr{A}}))$ is a $($topological or $C^0)$ Banach manifold built
  on the Banach space $C^{m+\mu}({\mathbb{S}}^{n-1})$ $($resp. $\dot{C}^{m+\mu}({\mathbb{S}}^{n-1}))$.}
\medskip

  {\em Proof}.\ \ We first note that for any $\Omega\in\mathfrak{D}^{m+\mu}({\mathbb{R}}^n)$ there exists corresponding
  local chart $(\mathcal{U},\varphi)\in\mathscr{A}$ such that $\Omega\in\mathcal{U}$. We omit the proof of this assertion
  here. Let $(\mathcal{U}_1,\varphi_1)$, $(\mathcal{U}_2,\varphi_2)$ be two local charts of
  $\mathfrak{D}^{m+\mu}({\mathbb{R}}^n)$ such that $\mathcal{U}_1\cap\mathcal{U}_2\neq\varnothing$, with base hypersurfaces
  $S_1$, $S_2$, respectively. Let $\bfn_i$ be the unit outward pointing normal field of $S_i$, $\mathcal{R}_i$ the
  neighborhood of $S_i$ and $\Pi_i:\mathcal{R}_i\to S_i$ the projection as defined above, $i=1,2$. Since $C^{m+\mu}(S_i)
  \approx C^{m+\mu}({\mathbb{S}}^{n-1})$, $i=1,2$, we only need to prove $\varphi_2\circ\varphi_1^{-1}\in
  C(\mathcal{O}_1,\mathcal{O}_2)$, where $\mathcal{O}_i=\varphi_i(\mathcal{U}_1\cap\mathcal{U}_2)\subseteq C^{m+\mu}(S_i)$,
  $i=1,2$. Let $\Omega\in\mathcal{U}_1\cap\mathcal{U}_2$. Then there exist $\rho_i\in\mathcal{O}_i$, $i=1,2$, such that
$$
  \partial\Omega=\{x+\rho_1(x)\bfn_1(x):x\in S_1\}=\{y+\rho_2(y)\bfn_2(y):y\in S_2\}.
$$
  By definition, $(\varphi_2\circ\varphi_1^{-1})(\rho_1)=\rho_2$, $\forall\rho_1\in\mathcal{O}_1$. The above equality implies
$$
  x+\rho_1(x)\bfn_1(x)=y+\rho_2(y)\bfn_2(y), \quad x\in S_1, \;\; y\in S_2,
$$
  i.e, for any $x\in S_1$ there exists a unique corresponding $y\in S_2$ such that the above equality holds, and vice
  versa. Since for given $\rho_1\in\mathcal{O}_1$ and $y\in S_2$, the point $x\in S_1$ is uniquely determined, we write
  $x=\chi(\rho_1,y)$. Substituting this expression into the above equality and computing inner products of both sides
  of it with $\bfn_2(y)$, we get
\begin{equation}
  \rho_2(y)=\langle\chi(\rho_1,y)-y,\bfn_2(y)\rangle+\rho_1(\chi(\rho_1,y))\langle\bfn_1(\chi(\rho_1,y)),\bfn_2(y)\rangle,
  \quad \forall y\in S_2.
\end{equation}
  Clearly, $x=\chi(\rho_1,y)$ is the implicit function defined by the equation
\begin{equation}
  y=\Pi_2(x+\rho_1(x)\bfn_1(x)).
\end{equation}
  This implicit function equation is regular, i.e., the derivative $\partial y/\partial x: T_x(S_1)\to T_y(S_2)$ has
  bounded inverse because $(\partial y/\partial x)^{-1}=\partial x/\partial y$ is the derivative of the function
  $y\mapsto x=\Pi_1(y+\rho_2(y)\bfn_2(y))$ which belongs to $C^{m+\mu}(S_2,S_1)$. Applying Lemma 2.1 we have
$$
  [(\rho_1,x)\mapsto\Pi_2(x+\rho_1(x)\bfn_1(x))]\in C^{m+\mu}(\mathcal{O}_1\times S_1,S_2).
$$
  It follows by the implicit function theorem that
\begin{equation}
  \chi\in C^{m+\mu}(\mathcal{O}_1\times S_2,S_1).
\end{equation}
  The above relation implies that
\begin{equation}
  [\rho_1\mapsto\chi(\rho_1,\cdot)]\in C^k(\mathcal{O}_1,C^{m-k+\mu}(S_2,S_1)), \quad k=0,1,\cdots,m.
\end{equation}
  Besides, by using (3.12) and applying Lemma 2.2, we have
\begin{equation}
  [(\rho_1,y)\mapsto\rho_1(\chi(\rho_1,y))]\in C^{m+\mu}(\mathcal{O}_1\times S_2,\mathbb{R}),
\end{equation}
  which further implies that
\begin{equation}
  [\rho_1\mapsto\rho_1(\chi(\rho_1,\cdot))]\in C^k(\mathcal{O}_1,C^{m-k+\mu}(S_2,\mathbb{R})), \quad k=0,1,\cdots,m.
\end{equation}
  Applying (3.13) and (3.15) for the case $k=0$, we obtain that the mapping $\rho_1\mapsto\rho_2$ belongs to
  $C(\mathcal{O}_1,\mathcal{O}_2)$, i.e., $\varphi_2\circ\varphi_1^{-1}\in C(\mathcal{O}_1,\mathcal{O}_2)$, as desired.
  Similarly, if $(\dot{\mathcal{U}}_1,\varphi_1)$, $(\dot{\mathcal{U}}_2,\varphi_2)$ be two local charts of
  $\dot{\mathfrak{D}}^{m+\mu}({\mathbb{R}}^n)$ such that $\dot{\mathcal{U}}_1\cap\dot{\mathcal{U}}_2\neq\varnothing$,
  then $\varphi_2\circ\varphi_1^{-1}\in C(\dot{\mathcal{O}}_1,\dot{\mathcal{O}}_2)$, where $\dot{\mathcal{O}}_i=
  \varphi_i(\dot{\mathcal{U}}_1\cap\dot{\mathcal{U}}_2)\subseteq\dot{C}^{m+\mu}(S_i)$, $i=1,2$. This proves the theorem.
  $\quad\Box$
\medskip

  Note that $\mathfrak{D}^{m+\mu}({\mathbb{R}}^n)$ (resp. $\dot{\mathfrak{D}}^{m+\mu}({\mathbb{R}}^n)$) does not possess
  a differentiable structure, because $\varphi_2\circ\varphi_1^{-1}$ is generally not differentiable. However, it is still
  possible to define differentiable points in $\mathfrak{D}^{m+\mu}({\mathbb{R}}^n)$ (resp.
  $\dot{\mathfrak{D}}^{m+\mu}({\mathbb{R}}^n)$) and tangent spaces at such points. These concepts will enable us to define
  differentiable curves in $\mathfrak{D}^{m+\mu}({\mathbb{R}}^n)$ (resp. $\dot{\mathfrak{D}}^{m+\mu}({\mathbb{R}}^n)$) and
  their tangent fields.
\medskip

  {\bf Lemma 3.5}\ \ {\em Let $(\mathcal{U}_1,\varphi_1)$, $(\mathcal{U}_2,\varphi_2)$ be two local charts of $\mathfrak{D}^{m+\mu}({\mathbb{R}}^n)$
  such that $\mathcal{U}_1\cap\mathcal{U}_2\neq\varnothing$.  Let $\mathcal{O}_i=\varphi_i(\mathcal{U}_1\cap\mathcal{U}_2)\subseteq C^{m+\mu}(S_i)$,
  $i=1,2$. Then for any $\Omega\in\mathcal{U}_1\cap\mathcal{U}_2\cap\mathfrak{D}^{m+\mu+1}({\mathbb{R}}^n)$, $\varphi_2\circ\varphi_1^{-1}$ is
  differentiable at $\rho=\varphi_1(\Omega)$: $(\varphi_2\circ\varphi_1^{-1})'(\rho)\in L(C^{m+\mu}(S_1),C^{m+\mu}(S_2))$, and
$$
  [\rho\mapsto (\varphi_2\circ\varphi_1^{-1})'(\rho)]\in C(\mathcal{O}_1\cap C^{m+\mu+1}(S_1),L(C^{m+\mu}(S_1),C^{m+\mu}(S_2))).
$$
  A similar result holds for $\dot{\mathfrak{D}}^{m+\mu}({\mathbb{R}}^n)$.}
\medskip

  {\em Proof}.\ \ Let $\sigma:\mathcal{O}_1\times S_2\to{\mathbb{R}}$ be the mapping given by the right-hand side of (3.3), i.e.,
$$
  \sigma(\rho,y)=\langle\chi(\rho,y)-y,\bfn_2(y)\rangle+\rho(\chi(\rho,y))\langle\bfn_1(\chi(\rho,y)),\bfn_2(y)\rangle, \quad
   \rho\in\mathcal{O}_1, \;\; y\in S_2,
$$
  where $\chi:\mathcal{O}_1\times S_2\to S_1$ is as before. Then $(\varphi_2\circ\varphi_1^{-1})(\rho)=\sigma(\rho,\cdot)$.
  From the above expression and the fact $\chi\in C^{m+\mu}(\mathcal{O}_1\times S_2,S_1)$ we see $\sigma\in C^{m+\mu}(\mathcal{O}_1\times S_2)$,
  by Lemma 2.2. Moreover, for any $\zeta\in C^{m+\mu}(S_1)$,
\begin{equation}
\begin{array}{rl}
   \partial_{\rho}\sigma(\rho,y)\zeta=&\langle[\partial_{\rho}\chi(\rho,y)\zeta](y),\bfn_2(y)\rangle+\{[\rho'(\chi(\rho,y))
   \partial_{\rho}\chi(\rho,y)\zeta](y)+\zeta(\chi(\rho,y))\}\langle\bfn_1(\chi(\rho,y)),\bfn_2(y)\rangle
\\ [0.1cm]
   &+\rho(\chi(\rho,y))\langle[\bfn_1'(\chi(\rho,y))\partial_{\rho}\chi(\rho,y)\zeta](y),\bfn_2(y)\rangle, \quad y\in S_2, \;\;
   \rho\in\mathcal{O}_1.
\end{array}
\end{equation}
  In what follows we prove:
\begin{equation}
  [\rho\mapsto[y\mapsto\partial_{\rho}\sigma(\rho,y)]]\in C(\mathcal{O}_1\cap C^{m+\mu+1}(S_1),C^{m+\mu}(S_2,L(C^{m+\mu}(S_1),{\mathbb{R}}))).
\end{equation}
  We first prove:
\begin{equation}
  [\rho\mapsto[y\mapsto\partial_{\rho}\chi(\rho,y)]]\in C(\mathcal{O}_1\cap C^{m+\mu+1}(S_1),C^{m+\mu}(S_2,L(C^{m+\mu}(S_1),{\mathbb{R}}^n))).
\end{equation}
  Indeed, letting $G:\mathcal{O}_1\times S_1\to S_2$ be the mapping given by the right-hand side of (3.11), we have
$$
  y=G(\rho,\chi(\rho,y)), \quad y\in S_2,\;\; \rho\in\mathcal{O}_1.
$$
  Differentiating both sides of the above equation in $\rho$, we get
\begin{equation}
  \partial_{\rho}\chi(\rho,y)=-[\partial_xG(\rho,\chi(\rho,y))]^{-1}\partial_{\rho}G(\rho,\chi(\rho,y)).
\end{equation}
  We have
\begin{equation}
  \partial_xG(\rho,x)\xi=\Pi_2'(x+\rho(x)\bfn_1(x))\{\xi+[\rho'(x)\xi]\bfn_1(x)+\rho(x)\bfn_1'(x)\xi\}, \quad \forall\xi\in T_x(S_1),
\end{equation}
\begin{equation}
  \partial_{\rho}G(\rho,x)\zeta=\Pi_2'(x+\rho(x)\bfn_1(x))[\zeta(x)\bfn_1(x)], \quad \forall\zeta\in C^{m+\mu}(S_1).
\end{equation}
  From (3.20) it is easy to see
$$
  [\rho\mapsto[(x,\xi)\mapsto(G(\rho,x),\partial_xG(\rho,x)\xi)]]\in C(\mathcal{O}_1\cap C^{m+\mu+1}(S_1),C^{m+\mu}(T(S_1),T(S_2))),
$$
  which combined with the fact $\chi\in C^{m+\mu}(\mathcal{O}_1\times S_2,S_1)$ implies
\begin{equation}
   [\rho\mapsto[(y,\eta)\mapsto(\chi(\rho,y),[\partial_xG(\rho,\chi(\rho,y))]^{-1}\eta)]]\in
   C(\mathcal{O}_1\cap C^{m+\mu+1}(S_1),C^{m+\mu}(T(S_2),T(S_1))).
\end{equation}
  From (3.21) it is also easy to see
\begin{equation}
   [\rho\mapsto[(x,\zeta)\mapsto(G(\rho,x),\partial_{\rho}G(\rho,x)\zeta)]]\in C(\mathcal{O}_1,C^{m+\mu}(S_1\times C^{m+\mu}(S_1),T(S_2))).
\end{equation}
  Combining (3.19), (3.22), (3.23) and using the fact $\chi\in C^{m+\mu}(\mathcal{O}_1\times S_2,S_1)\subseteq C(\mathcal{O}_1,C^{m+\mu}(S_2,S_1))$
  we see that (3.18) follows. From (3.16), (3.18) and the fact $\chi\in C^{m+\mu}(\mathcal{O}_1\times S_2,S_1)\subseteq C(\mathcal{O}_1,
  C^{m+\mu}(S_2,S_1))$ we obtain (3.17). Hence
$$
  [\rho\mapsto(\varphi_2\circ\varphi_1^{-1})'(\rho)]=[\rho\mapsto[\zeta\mapsto\partial_{\rho}\sigma(\rho,\cdot)\zeta]]\in
  C(\mathcal{O}_1\cap C^{m+\mu+1}(S_1),L(C^{m+\mu}(S_1),C^{m+\mu}(S_2))).
$$
  This proves the desired assertion. $\quad\Box$
\medskip

  We recall the following concept given in \cite{Cui4}:
\medskip

  {\bf Definition 3.6}\ {\em Let $\mathfrak{M}$ and $\mathfrak{M}_0$ be two Banach manifolds built on Banach spaces $X$ and $X_0$,
  respectively. Let $\mathscr{A}$ be a family of local charts of $\mathfrak{M}$. We say $\mathfrak{M}_0$ is a \mbox{\small\boldmath
  $C^1$-$embedded\;Banach\; submanifold$} of $\mathfrak{M}$ with respect to $\mathscr{A}$ if the following conditions are satisfied:
\begin{enumerate}
\item[]$(1)$\ \ $X_0$ is a densely embedded Banach subspace of $X$.\vspace*{-0.2cm}
\item[]$(2)$\ \ $\mathfrak{M}_0$ is an embedded topological subspace of $\mathfrak{M}$, i.e., $\mathfrak{M}_0\subseteq
  \mathfrak{M}$, and for any open subset $U$ of $\mathfrak{M}$, $U\cap\mathfrak{M}_0$ is an open subset of $\mathfrak{M}_0$.
\vspace*{-0.2cm}
\item[]$(3)$\ \ For any $\eta\in\mathfrak{M}$ there exists a local chart $(\mathcal{U},\varphi)\in\mathscr{A}$ such
  that $\eta\in\mathcal{U}$.
\item[]$(4)$\ \ For any $\eta\in\mathfrak{M}_0$ and any local chart $(\mathcal{U},\varphi)\in\mathscr{A}$ such that
  $\eta\in\mathcal{U}$, by letting $\mathcal{U}_0=\mathcal{U}\cap\mathfrak{M}_0$, $(\mathcal{U}_0,\varphi|_{\mathcal{U}_0})$
  is a local chart of $\mathfrak{M}_0$ at $\eta$.  \vspace*{-0.2cm}
\item[]$(5)$\ \ For any $\eta\in\mathfrak{M}_0$ and any $(\mathcal{U},\varphi),(\mathcal{V},\psi)\in\mathscr{A}$ such that
  $\eta\in\mathcal{U}\cap\mathcal{V}$, by letting $\mathcal{U}_0=\mathcal{U}\cap\mathfrak{M}_0$ and $\mathcal{V}_0=\mathcal{V}\cap
  \mathfrak{M}_0$, the following relations hold:
$$
  \psi\circ\varphi^{-1}\in\mathfrak{C}^1(\varphi(\mathcal{U}_0\cap\mathcal{V}_0);X,X) \quad \mbox{and} \quad
  \varphi\circ\psi^{-1}\in\mathfrak{C}^1(\psi(\mathcal{U}_0\cap\mathcal{V}_0);X,X).
$$
\end{enumerate}
  We call any local chart $(\mathcal{U},\varphi)\in\mathscr{A}$ such that $\eta\in\mathcal{U}$ a \mbox{\small\boldmath
  $(C^1,\mathfrak{M}_0)$-$regular\;local\; chart$} of $\mathfrak{M}$ at $\eta$, and call the family $\mathscr{A}$ a
  \mbox{\small\boldmath $(C^1,\mathfrak{M}_0)$-$regular$} or simply a \mbox{\small\boldmath $\mathfrak{M}_0$-$regular\;
  local\; chart\;system$} of $\mathfrak{M}$.}
\medskip

  As an immediate consequence of Lemma 3.3, we have the following result:
\medskip

  {\bf Theorem 3.7}\ {\em Let $m$ be a positive integer and $0\leqslant\mu\leqslant 1$. Let $\mathscr{A}$ be the set of all
  local charts of $\dot{\mathfrak{D}}^{m+\mu}({\mathbb{R}}^n)$ as considered in Theorem 3.2. With respect to this local
  chart system, $\dot{\mathfrak{D}}^{m+\mu+1}({\mathbb{R}}^n)$ is a $C^1$-embedded Banach submanifold of
  $\dot{\mathfrak{D}}^{m+\mu}({\mathbb{R}}^n)$.} $\quad\Box$
\medskip

  We also recall the following concepts given in \cite{Cui4}:
\medskip

  {\bf Definition 3.8}\  {\em Let $\mathfrak{M}$ be a Banach manifold, $\mathscr{A}$ a family of local charts of $\mathfrak{M}$,
  and $\mathfrak{M}_0$ a $C^1$-embedded Banach submanifold of $\mathfrak{M}$ with respect to $\mathscr{A}$. We have the following
  notions:

  $(1)$\ We say $(\mathfrak{M},\mathfrak{M}_0,\mathscr{A})$ is \mbox{\small\boldmath $inward\;spreadable$} if there exists a
  Banach manifold $\mathfrak{M}_1\subseteq\mathfrak{M}_0$ such that $\mathfrak{M}_1$ is a $C^1$-embedded Banach submanifold
  of $\mathfrak{M}_0$ with respect to the restriction of $\mathscr{A}$ to $\mathfrak{M}_0$. In this case, a local chart in
  $\mathscr{A}$ is called a $(\mathfrak{M}_0,\mathfrak{M}_1)$-regular local chart.

  $(2)$\ We say $(\mathfrak{M},\mathfrak{M}_0,\mathscr{A})$ is \mbox{\small\boldmath $outward\;spreadable$} if there exists a
  Banach manifold $\widetilde{\mathfrak{M}}\supseteq\mathfrak{M}$ and a family $\widetilde{\mathscr{A}}$ of local charts of
  $\widetilde{\mathfrak{M}}$, such that $\mathscr{A}$ is the restriction of $\widetilde{\mathscr{A}}$ to $\mathfrak{M}$ and
  $\mathfrak{M}$ is a $C^1$-embedded Banach submanifold of $\widetilde{\mathfrak{M}}$ with respect to $\widetilde{\mathscr{A}}$.

  $(3)$\ If $(\mathfrak{M},\mathfrak{M}_0,\mathscr{A})$ is both inward spreadable and outward spreadable then we call the pair
  $(\mathfrak{M},\mathscr{A})$ a \mbox{\small\boldmath $quasi$-$di\!f\!ferentiable$} Banach manifold with a \mbox{\small\boldmath
  $C^1$-$kernel$} $\mathfrak{M}_0$, or simply call it a quasi-differentiable Banach manifold without mentioning the
  $C^1$-kernel $\mathfrak{M}_0$. $\mathfrak{M}_1$ is called an \mbox{\small\boldmath $inner$ $C^1$-$kernel$} of
  $(\mathfrak{M},\mathscr{A})$, and $(\widetilde{\mathfrak{M}},\widetilde{\mathscr{A}})$ a \mbox{\small\boldmath $C^1$-$shell$}
  of $(\mathfrak{M},\mathscr{A})$.}
\medskip

  The following result is an immediate consequence of Theorem 3.7:
\medskip

  {\bf Theorem 3.9}\ {\em Let $m$ be a positive integer not less than $2$, and $0\leqslant\mu\leqslant 1$. Let $\mathscr{A}$
  be the set of all local charts of $\dot{\mathfrak{D}}^{m+\mu}({\mathbb{R}}^n)$ as considered in Theorem 3.2. Then
  $(\dot{\mathfrak{D}}^{m+\mu}({\mathbb{R}}^n),\mathscr{A})$ is a quasi-differentiable Banach manifold with a $C^1$-kernel
  $\dot{\mathfrak{D}}^{m+1+\mu}({\mathbb{R}}^n)$, an inner $C^1$-kernel $\dot{\mathfrak{D}}^{m+2+\mu}({\mathbb{R}}^n)$, and
  a shell $\dot{\mathfrak{D}}^{m-1+\mu}({\mathbb{R}}^n)$.} $\quad\Box$
\medskip

  From \cite{Cui4} we know that if $(\mathfrak{M},\mathscr{A})$ is a quasi-differentiable Banach manifold and $\mathfrak{M}_0$
  is a $C^1$-kernel of it, then for any $x\in\mathfrak{M}_0$ the tangent space $\mathcal{T}_x(\mathfrak{M})$ makes sense,
  which is a Banach space isomorphic to the base Banach space $X$ of $\mathfrak{M}$. Let us briefly recall this notion for
  the Banach manifold $\dot{\mathfrak{D}}^{m+\mu}({\mathbb{R}}^n)$.
\medskip

  {\bf Definition 3.10}\ \ {\em Let $m$ be a positive integer not less than $2$, and $0\leqslant\mu\leqslant 1$. Let
  $\Omega\in\dot{\mathfrak{D}}^{m+\mu+1}({\mathbb{R}}^n)$, and regard it as a point in $\dot{\mathfrak{D}}^{m+\mu}({\mathbb{R}}^n)$.
  We have the following concepts:

  $(1)$\ \ For a function $F:\mathcal{O}\to{\mathbb{R}}$ defined in a neighborhood $\mathcal{O}\subseteq
  \dot{\mathfrak{D}}^{m+\mu}({\mathbb{R}}^n)$ of $\Omega$, we say $F$ is \mbox{\small\boldmath $fully$
  $strongly$ $continuously$ $di\!f\!ferentiable$} at $\Omega$ if there exists a neighborhood $\mathcal{O}'\subseteq\mathcal{O}$
  of $\Omega$ in $\dot{\mathfrak{D}}^{m+\mu}({\mathbb{R}}^n)$ such that for any regular local chart $(\mathcal{U},\varphi)$ of
  $\dot{\mathfrak{D}}^{m+\mu}({\mathbb{R}}^n)$ at $\Omega$, the function $F\circ\varphi^{-1}:\varphi(\mathcal{O}\cap\mathcal{U})
  \to{\mathbb{R}}$ is continuously differentiable in $O=\varphi(\mathcal{O}'\cap\mathcal{U})$ in the topology of
  $X=\dot{C}^{m+\mu}(S)$ and $[u\mapsto(F\circ\varphi^{-1})'(u)]\in C(O,X^*)=C(O,L(X,\mathbb{R}))$, where $S$ is the base
  hypersurface of the local chart $(\mathcal{U},\varphi)$ and $O$ uses the topology of $X$. We denote by
  $\dot{\mathscr{D}}^{1\!s}_{\Omega}$ the set of all real-valued functions $F$ defined in a neighborhood of $\Omega$ in
  $\dot{\mathfrak{D}}^{m+\mu}({\mathbb{R}}^n)$ which are fully strongly continuously differentiable at $\Omega$.

  $(2)$\ \ Let $f:(-\varepsilon,\varepsilon)\to\dot{\mathfrak{D}}^{m+\mu}({\mathbb{R}}^n)$ $(\varepsilon>0)$ be a curve in
  $\dot{\mathfrak{D}}^{m+\mu}({\mathbb{R}}^n)$ passing $\Omega$, i.e., $f(0)=\Omega$. We say $f(t)$ is \mbox{\small\boldmath
  $differentiable$} at $t=0$ if there exists a regular local chart $(\mathcal{U},\varphi)$ of
  $\dot{\mathfrak{D}}^{m+\mu}({\mathbb{R}}^n)$ at $\Omega$, such that the function $t\mapsto\varphi(f(t))$ is differentiable
  at $t=0$ in the topology of $X=\dot{C}^{m+\mu}(S)$. Moreover, we define the \mbox{\small\boldmath $tangent$ $vector$} of
  this curve at $\Omega$, or the \mbox{\small\boldmath $derivative$} $f'(0)$ of $f(t)$ at $t=0$, to be the mapping
  $f'(0):\dot{\mathscr{D}}^{1\!s}_{\Omega}\to{\mathbb{R}}$ defined by
$$
  f'(0)F=(F\circ f)'(0), \quad \forall F\in\dot{\mathscr{D}}^{1\!s}_{\Omega}.
$$

  $(3)$\ \ We denote
$$
  T_{\Omega}(\dot{\mathfrak{D}}^{m+\mu}({\mathbb{R}}^n))=\{f'(0):f:(-\varepsilon,\varepsilon)
  \to\dot{\mathfrak{D}}^{m+\mu}({\mathbb{R}}^n),f(0)=\Omega,
  \mbox{$f(t)$ is differentiable at $t=0$}\},
$$
  and call it the \mbox{\small\boldmath $tangent$ $space$} of $\dot{\mathfrak{D}}^{m+\mu}({\mathbb{R}}^n)$ at $\Omega$.}
\medskip

  From the discussion made in \cite{Cui4} we know that the above concepts make sense, and
  $T_{\Omega}(\dot{\mathfrak{D}}^{m+\mu}({\mathbb{R}}^n))$ is a Banach space isomorphic to $X=\dot{C}^{m+\mu}(S)$.
\medskip

\section{Standard local chart of $\mathfrak{D}^{m+\mu}({\mathbb{R}}^n)$}
\setcounter{equation}{0}

\hskip 2em
  For $\Omega\in{\mathfrak{D}}^{m+\mu+1}({\mathbb{R}}^n)$, the tangent space $T_{\Omega}({\mathfrak{D}}^{m+\mu}({\mathbb{R}}^n))$
  can be expressed in a different form which is very useful from the viewpoint of application. To get that expression, let us first
  consider a local chart of ${\mathfrak{D}}^{m+\mu}({\mathbb{R}}^n)$ at a point $\Omega\in{\mathfrak{D}}^{m+\mu+1}({\mathbb{R}}^n)$
  which is different from those introduced before. Indeed, since $\Omega\in{\mathfrak{D}}^{m+\mu+1}({\mathbb{R}}^n)$ implies that
  $\Gamma:=\partial\Omega$ is a $C^{m+\mu+1}$-hypersurface and its normal field $\bfnu$ is of $C^{m+\mu}$-class: $\bfnu\in
  C^{m+\mu}(\Gamma,{\mathbb{R}}^n)$, by using $\Gamma$ as a base hypersurface and repeating the argument before, we get a
  (irregular) local chart of ${\mathfrak{D}}^{m+\mu}({\mathbb{R}}^n)$ at $\Omega$, which we denote as
  $(\mathcal{U}_{\Omega},\varphi_{\Omega})$.
\medskip

  {\bf Definition 4.1}\ \ {\em We call $(\mathcal{U}_{\Omega},\varphi_{\Omega})$ the  \mbox{\small\boldmath $standard$ $local$
  $chart$} of $\mathfrak{D}^{m+\mu}({\mathbb{R}}^n)$ at $\Omega$.}
\medskip

  {\bf Lemma 4.2}\ \ {\em Let $\Omega\in{\mathfrak{D}}^{m+\mu+1}({\mathbb{R}}^n)$ and $\Gamma=\partial\Omega$. For any
  regular local chart $(\mathcal{U},\varphi)$ of ${\mathfrak{D}}^{m+\mu}({\mathbb{R}}^n)$ at $\Omega$ with a smooth base
  hypersurface $S$, the coordinate transformation mappings $\varphi_{\Omega}\circ\varphi^{-1}$ and
  $\varphi\circ\varphi_{\Omega}^{-1}$ satisfy the following properties:
\begin{equation}
  \varphi_{\Omega}\circ\varphi^{-1}\in C(\varphi(\mathcal{U}\cap\mathcal{U}_{\Omega}),\varphi_{\Omega}(\mathcal{U}\cap\mathcal{U}_{\Omega})),
   \quad
  \varphi\circ\varphi_{\Omega}^{-1}\in C(\varphi_{\Omega}(\mathcal{U}\cap\mathcal{U}_{\Omega}),\varphi(\mathcal{U})\cap\mathcal{U}_{\Omega})).
\end{equation}
  If further $\Omega\in{\mathfrak{D}}^{m+\mu+2}({\mathbb{R}}^n)$ then $\varphi_{\Omega}\circ\varphi^{-1}$ and
  $\varphi\circ\varphi_{\Omega}^{-1}$ are differentiable at the points in $\varphi(\mathcal{U}_1\cap\mathcal{U}_{\Omega1})$
  and $\varphi_{\Omega}(\mathcal{U}_1\cap\mathcal{U}_{\Omega1})$ respectively, where $\mathcal{U}_1=\mathcal{U}\cap
  {\mathfrak{D}}^{m+\mu+1}({\mathbb{R}}^n)$, $\mathcal{U}_{\Omega1}=\mathcal{U}_{\Omega}\cap{\mathfrak{D}}^{m+\mu+1}({\mathbb{R}}^n)$,
  and
\begin{equation}
\left\{
\begin{array}{l}
   (\varphi_{\Omega}\circ\varphi^{-1})'\in C(\varphi(\mathcal{U}_1\cap\mathcal{U}_{\Omega1}),L({C}^{m+\mu}(S),{C}^{m+\mu}(\Gamma))),
    \\
   (\varphi\circ\varphi_{\Omega}^{-1})'\in C(\varphi_{\Omega}(\mathcal{U}_1\cap\mathcal{U}_{\Omega1}),L({C}^{m+\mu}(\Gamma),{C}^{m+\mu}(S))).
\end{array}
\right.
\end{equation}
  Moreover $($recall that $\varphi_{\Omega}(\Omega)=0)$,}
\begin{equation}
   (\varphi_{\Omega}\circ\varphi^{-1})'(\varphi(\Omega))\in L({C}^{m+\mu+1}(S),{C}^{m+\mu+1}(\Gamma)), \quad
   (\varphi\circ\varphi_{\Omega}^{-1})'(0)\in L({C}^{m+\mu+1}(\Gamma),{C}^{m+\mu+1}(S)).
\end{equation}

  {\em Proof}.\ \ Let $\bfn$ be the unit normal field of $S$, outward pointing with respect to the domain $Q$ enclosed by $S$.
  Let $\bfnu$ be as above, i.e., it is the unit normal field of $\Gamma$, outward pointing with respect to $\Omega$. Let
  $\mathcal{R}$ be a neighborhood of $S$ and $\Pi:\mathcal{R}\to S$ the projection, both defined as before. Let $\mathcal{R}_{\Omega}$
  be a neighborhood of $\Gamma$ similar to $\mathcal{R}$ and $\Pi_{\Omega}:\mathcal{R}_{\Omega}\to\Gamma$ the projection
  similar to $\Pi$. Let $x$, $y$, $\rho$ and $\eta$ denote variables in $S$, $\Gamma$, ${C}^{m+\mu}(S)$ and
  ${C}^{m+\mu}(\Gamma)$, respectively. Let $\Delta\in\mathcal{U}\cap\mathcal{U}_{\Omega}$. Then there exist $\rho\in
  \varphi(\mathcal{U}\cap\mathcal{U}_{\Omega})$ and $\eta\in\varphi_{\Omega}(\mathcal{U}\cap\mathcal{U}_{\Omega})$ such that
$$
  \partial\Delta=\{x+\rho(x)\bfn(x):x\in S\}=\{y+\eta(y)\bfnu(y):y\in\Gamma\}.
$$
  Then $(\varphi_{\Omega}\circ\varphi^{-1})(\rho)=\eta$, $\forall\rho\in\varphi(\mathcal{U}\cap\mathcal{U}_{\Omega})$,
  and $(\varphi\circ\varphi_{\Omega}^{-1})(\eta)=\rho$, $\forall\eta\in\varphi_{\Omega}(\mathcal{U}\cap\mathcal{U}_{\Omega})$.
  Let $x=\chi(\rho,y)$ and $y=\tau(\eta,x)$ be the implicit functions defined by the equations
\begin{equation}
  y=\Pi_{\Omega}(x+\rho(x)\bfn(x)) \quad \mbox{and} \quad x=\Pi(y+\eta(y)\bfnu(y)),
\end{equation}
  respectively. From the proof of Theorem 3.4 we see that
\begin{equation}
  \eta(y)=\langle\chi(\rho,y)-y,\bfnu(y)\rangle+\rho(\chi(\rho,y))\langle\bfn(\chi(\rho,y)),\bfnu(y)\rangle,
  \quad y\in\Gamma,
\end{equation}
\begin{equation}
  \rho(x)=\langle\tau(\eta,x)-x,\bfn(x)\rangle+\eta(\tau(\eta,x))\langle\bfnu(\tau(\eta,x)),\bfn(x)\rangle,
  \quad x\in S.
\end{equation}
  We have seen that $\bfn\in C^{\infty}(S,{\mathbb{R}}^n)$, $\Pi\in C^{\infty}(\mathcal{R},S)$ and
  $\bfnu\in{C}^{m+\mu}(\Gamma,{\mathbb{R}}^n)$. Moreover, it is easy to see that $\Pi_{\Omega}\in
  {C}^{m+\mu}(\mathcal{R}_{\Omega},\Gamma)$. From these facts and some similar argument as in the proof of
  Theorem 3.4 we see that $\chi\in{C}^{m+\mu}(\varphi(\mathcal{U}\cap\mathcal{U}_{\Omega})\times\Gamma,S)$
  and $\tau\in{C}^{m+\mu}(\varphi_{\Omega}(\mathcal{U}\cap\mathcal{U}_{\Omega})\times S,\Gamma)$. Consequently,
  the two relations in (4.1) follow.

  Next assume that $\Omega\in{\mathfrak{D}}^{m+\mu+2}({\mathbb{R}}^n)$. Then $\bfnu\in{C}^{m+\mu+1}(\Gamma,{\mathbb{R}}^n)$.
  Given $\rho\in\varphi(\mathcal{U}_1\cap\mathcal{U}_{\Omega1})$, $\eta\in\varphi_{\Omega}(\mathcal{U}_1\cap\mathcal{U}_{\Omega1})$,
  $\xi\in{C}^{m+\mu}(S)$ and $\zeta\in{C}^{m+\mu}(\Gamma)$, we denote $u=(\varphi_{\Omega}\circ\varphi^{-1})'(\rho)\xi$
  and $v=(\varphi\circ\varphi_{\Omega}^{-1})'(\eta)\zeta$. Then similarly as in the proof of Lemma 3.5 we have
\begin{equation}
\begin{array}{rl}
   u(y)=&\langle\partial_{\rho}\chi(\rho,y)\xi,\bfnu(y)\rangle+[\rho'(\chi(\rho,y))
   \partial_{\rho}\chi(\rho,y)\xi+\xi(\chi(\rho,y))]\langle\bfn(\chi(\rho,y)),\bfnu(y)\rangle
\\ [0.1cm]
   &+\rho(\chi(\rho,y))\langle\bfn'(\chi(\rho,y))\partial_{\rho}\chi(\rho,y)\xi,\bfnu(y)\rangle, \quad y\in\Gamma,
\end{array}
\end{equation}
\begin{equation}
\begin{array}{rl}
   v(x)=&\langle\partial_{\eta}\tau(\eta,x)\zeta,\bfn(x)\rangle+[\eta'(\tau(\eta,x))
   \partial_{\eta}\tau(\eta,x)\zeta+\zeta(\tau(\eta,x))]\langle\bfnu(\tau(\eta,x)),\bfn(x)\rangle
\\ [0.1cm]
   &+\eta(\tau(\eta,x))\langle\bfnu'(\tau(\eta,x))\partial_{\eta}\tau(\eta,x)\zeta,\bfn(x)\rangle, \quad x\in S.
\end{array}
\end{equation}
  From these relations we easily see that the two relations in (4.2) hold. In particular, if $\rho=\varphi(\Omega)$
  then $\rho\in{C}^{m+\mu+2}(S)$. Since $\bfn\in C^{\infty}(S)$, and the condition $\Omega\in{\mathfrak{D}}^{m+\mu+2}({\mathbb{R}}^n)$
  implies that $\bfnu\in{C}^{m+\mu+1}(\Gamma,{\mathbb{R}}^n)$ and $\chi\in{C}^{m+\mu+1}(\varphi(\mathcal{U}_1\cap
  \mathcal{U}_{\Omega1})\times\Gamma,S)$, from (4.7) we easily see that the first relation in (4.3) holds. Moreover,
  if $\eta=0$ then (4.8) becomes
\begin{equation}
  v(x)=\langle\partial_{\eta}\tau(0,x)\zeta,\bfn(x)\rangle+\zeta(\tau(0,x))
  \langle\bfnu(\tau(0,x)),\bfn(x)\rangle, \quad x\in S.
\end{equation}
  Let $G:\varphi_{\Omega}(\mathcal{U}_1\cap\mathcal{U}_{\Omega1})\times\Gamma\to S$ be the map defined by
  $G(\eta,y)=\Pi(y+\eta(y)\bfnu(y))$, $\eta\in\varphi_{\Omega}(\mathcal{U}_1\cap\mathcal{U}_{\Omega1})$, $y\in\Gamma$.
  We have
$$
  x=G(\eta,\tau(\eta,x)), \quad x\in S,\;\; \eta\in\varphi_{\Omega}(\mathcal{U}_1\cap\mathcal{U}_{\Omega1}).
$$
  Differentiating both sides of the above equation in $\eta$, we get
$$
  \partial_{\eta}\tau(\eta,x)=-[\partial_yG(\eta,\tau(\eta,x))]^{-1}\partial_{\eta}G(\eta,\tau(\eta,x)).
$$
  In particular, since $\tau(0,x)=\rho_0(x)$, where $\rho_0=\varphi(\Omega)$, we get
\begin{equation}
  \partial_{\eta}\tau(0,x)=-[\partial_yG(0,\rho_0(x))]^{-1}\partial_{\eta}G(0,\rho_0(x)).
\end{equation}
  We have
$$
  \partial_yG(\eta,y)z=\Pi'(y+\eta(y)\bfnu(y))\{z+[\eta'(y)z]\bfnu(y)+\eta(y)\bfnu'(y)z\},
  \quad \forall z\in T_y(\Gamma),
$$
$$
  \partial_{\eta}G(\eta,y)\zeta=\Pi'(y+\eta(y)\bfnu(y))[\zeta(y)\bfnu(y)], \quad \forall\zeta\in C^{m+\mu}(\Gamma).
$$
  In particular,
$$
  \partial_yG(0,y)z=\Pi'(y)z, \quad \forall z\in T_y(\Gamma); \quad
  \partial_{\eta}G(0,y)\zeta=\Pi'(y)[\zeta(y)\bfnu(y)], \quad \forall\zeta\in C^{m+\mu}(\Gamma).
$$
  Note that $\Pi'(y)$ is not injective, so that it does not have an inverse. Let $P:\Gamma\to L(\mathbb{R}^n)$ be
  the map defined by $P(y)z=\langle z,\bfnu(y)\rangle\bfnu(y)$, $y\in\Gamma$, $z\in\mathbb{R}^n$. Then
\begin{equation}
  \partial_yG(0,y)=A(y), \quad  \partial_{\eta}G(0,y)=B(y), \quad \forall y\in\Gamma,
\end{equation}
  where for every $y\in\Gamma$,
\begin{equation}
  A(y)z=\Pi'(y)[I-P(y)]z, \quad \forall z\in T_y(\Gamma); \quad
  B(y)\zeta=\Pi'(y)P(y)[\zeta(y)\bfnu(y)], \quad \forall\zeta\in C^{m+\mu}(\Gamma).
\end{equation}
  Clearly, when $S$ is sufficiently close to $\Gamma$ which we assume to be true, for every $y\in\Gamma$,
  $A(y):T_y(\Gamma)\to T_{\Pi(y)}(S)$ is invertible. By (4.10) and (4.11), it follows that $\partial_{\eta}\tau(0,x)
  =-A(\rho_0(x))^{-1}B(\rho_0(x))$, $\forall x\in S$. Hence, from (4.9) we get
\begin{equation}
  v(x)=-\langle A(\rho_0(x))^{-1}B(\rho_0(x))\zeta,\bfn(x)\rangle+\zeta(\rho_0(x))
  \langle\bfnu(\rho_0(x)),\bfn(x)\rangle, \quad x\in S.
\end{equation}
  From this expression and the definition (4.12) of the operators $A$ and $B$ we easily see that the second relation
  in (4.3) holds. This completes the proof of the lemma. $\quad\Box$
\medskip

  {\bf Corollary 4.3}\ \ {\em Let $\Omega\in{\mathfrak{D}}^{m+\mu+1}({\mathbb{R}}^n)$ and $\Gamma=\partial\Omega$,
  where $m\geqslant 2$. For any regular local chart $(\mathcal{U},\varphi)$ of ${\mathfrak{D}}^{m+\mu}({\mathbb{R}}^n)$
  at $\Omega$ with a smooth base hypersurface $S$, the coordinate transformation mappings $\varphi_{\Omega}\circ\varphi^{-1}$ and
  $\varphi\circ\varphi_{\Omega}^{-1}$ are differentiable at the points in $\varphi(\mathcal{U}\cap\mathcal{U}_{\Omega})$
  and $\varphi_{\Omega}(\mathcal{U}\cap\mathcal{U}_{\Omega})$ respectively, in the topologies of ${C}^{m+\mu-1}(S)$ and
  ${C}^{m+\mu-1}(\Gamma)$, and
\begin{equation}
\left\{
\begin{array}{l}
   (\varphi_{\Omega}\circ\varphi^{-1})'\in C(\varphi(\mathcal{U}\cap\mathcal{U}_{\Omega}),L({C}^{m+\mu-1}(S),{C}^{m+\mu-1}(\Gamma))),
    \\
   (\varphi\circ\varphi_{\Omega}^{-1})'\in C(\varphi_{\Omega}(\mathcal{U}\cap\mathcal{U}_{\Omega}),L({C}^{m+\mu-1}(\Gamma),{C}^{m+\mu-1}(S))).
\end{array}
\right.
\end{equation}
  Moreover $($recall that $\varphi_{\Omega}(\Omega)=0)$,}
\begin{equation}
   (\varphi_{\Omega}\circ\varphi^{-1})'(\varphi(\Omega))\in L({C}^{m+\mu}(S),{C}^{m+\mu}(\Gamma)), \quad
   (\varphi\circ\varphi_{\Omega}^{-1})'(0)\in L({C}^{m+\mu}(\Gamma),{C}^{m+\mu}(S)).
\end{equation}
  $\Box$
\medskip

  Since even if $\Omega\in{\mathfrak{D}}^{m+\mu+1}({\mathbb{R}}^n)$, the normal field $\bfnu$ of $\Gamma=\partial\Omega$
  is not differentiable in $C^{m+\mu}(\Gamma)$, in the standard local chart tangent vector and tangent space at $\Omega$
  cannot be defined similarly as in regular local chart. We give the following alternative definition:
\medskip

  {\bf Definition 4.4}\ \ {\em Let $m$ be a positive integer not less than $2$, and $0\leqslant\mu\leqslant 1$. Let
  $\Omega\in{\mathfrak{D}}^{m+\mu+1}({\mathbb{R}}^n)$, and regard it as a point in ${\mathfrak{D}}^{m+\mu}({\mathbb{R}}^n)$.
  We have the following concepts:

  $(1)$\ \ For a function $F:\mathcal{O}\to{\mathbb{R}}$ defined in a neighborhood $\mathcal{O}\subseteq
  {\mathfrak{D}}^{m+\mu}({\mathbb{R}}^n)$ of $\Omega$, we say $F$ is \mbox{\small\boldmath $very$ $strongly$
  $continuously$ $di\!f\!ferentiable$} at $\Omega$ if there exists a neighborhood $\mathcal{O}'\subseteq\mathcal{O}$
  of $\Omega$ in ${\mathfrak{D}}^{m+\mu}({\mathbb{R}}^n)$ such that the function $F\circ\varphi_{\Omega}^{-1}:
  \varphi_{\Omega}(\mathcal{O}\cap\mathcal{U})\to{\mathbb{R}}$ is continuously differentiable in $O=\varphi_{\Omega}(\mathcal{O}'
  \cap\mathcal{U})$ in the topology of $\tilde{X}={C}^{m-1+\mu}(\Gamma)$ and $[u\mapsto(F\circ\varphi_{\Omega}^{-1})'(u)]
  \in C(O,\tilde{X}^*)=C(O,L(\tilde{X},\mathbb{R}))$, where $O$ uses the topology of $X={C}^{m+\mu}(\Gamma)$. We denote by
  $\hat{\mathscr{D}}^{1\!s}_{\Omega}$ the set of all real-valued functions $F$ defined in a neighborhood of $\Omega$ in
  ${\mathfrak{D}}^{m+\mu}({\mathbb{R}}^n)$ which are very strongly continuously differentiable at $\Omega$.

  $(2)$\ \ Let $f:(-\varepsilon,\varepsilon)\to{\mathfrak{D}}^{m+\mu}({\mathbb{R}}^n)$ $(\varepsilon>0)$ be a curve in
  ${\mathfrak{D}}^{m+\mu}({\mathbb{R}}^n)$ passing $\Omega$, i.e., $f(0)=\Omega$. We say $f(t)$ is \mbox{\small\boldmath
  $weakly$ $differentiable$} at $t=0$ if the function $t\mapsto\varphi_{\Omega}(f(t))$ is differentiable at $t=0$ in the
  topology of $\tilde{X}={C}^{m-1+\mu}(\Gamma)$ and $(\varphi_{\Omega}\circ f)'(0)\in X={C}^{m+\mu}(\Gamma)$.
  Moreover, we define the \mbox{\small\boldmath $tangent$ $vector$} of this curve at $\Omega$, or the
  \mbox{\small\boldmath $derivative$} $f'(0)$ of $f(t)$ at $t=0$, to be the mapping $f'(0):{\mathscr{D}}^{1\!s}_{\Omega}
  \to{\mathbb{R}}$ defined by}
\begin{equation}
  f'(0)F=(F\circ f)'(0), \quad \forall F\in\hat{\mathscr{D}}^{1\!s}_{\Omega}.
\end{equation}

  {\em Remark}.\ \ Since
$$
  F\circ f=(F\circ\varphi_{\Omega}^{-1})\circ(\varphi_{\Omega}\circ f),
$$
  we see that the right-hand side of (4.16) makes sense.
\medskip

  {\bf Lemma 4.5}\ \ {\em We have the following assertions:

  $(1)$\ \ For a function $F:\mathcal{O}\to{\mathbb{R}}$ defined in a neighborhood $\mathcal{O}\subseteq
  {\mathfrak{D}}^{m+\mu}({\mathbb{R}}^n)$ of $\Omega$, if $F$ is very strongly continuously differentiable at $\Omega$
  then $F$ is also fully strongly continuously differentiable at $\Omega$.

  $(2)$\ \ For a curve $f:(-\varepsilon,\varepsilon)\to{\mathfrak{D}}^{m+\mu}({\mathbb{R}}^n)$ $(\varepsilon>0)$
  in ${\mathfrak{D}}^{m+\mu}({\mathbb{R}}^n)$ passing $\Omega$, i.e., $f(0)=\Omega$, if $f$ is differentiable at
  $\Omega$ then $f$ is also weakly differentiable at $\Omega$.

  $(3)$\ \ The following relation holds: }
\begin{equation}
  T_{\Omega}({\mathfrak{D}}^{m+\mu}({\mathbb{R}}^n))=\{f'(0):f:(-\varepsilon,\varepsilon)
  \to{\mathfrak{D}}^{m+\mu}({\mathbb{R}}^n),f(0)=\Omega,
  \mbox{$f(t)$ is weakly differentiable at $t=0$}\}.
\end{equation}

  {\em Proof}.\ \ (1)\ Let $({\mathcal{U}},\varphi)$ be a regular local chart of ${\mathfrak{D}}^{m+\mu}({\mathbb{R}}^n)$
  at $\Omega$. We have
$$
  F\circ\varphi^{-1}=(F\circ\varphi_{\Omega}^{-1})\circ(\varphi_{\Omega}\circ\varphi^{-1}).
$$
  By Corollary 4.3, $\varphi_{\Omega}\circ\varphi^{-1}$ is continuously differentiable as a map from
  $\varphi({\mathcal{U}}\cap\mathcal{U}_{\Omega})\subseteq {C}^{m-1+\mu}(S)$ to ${C}^{m-1+\mu}(\Gamma)$, and
  $(\varphi_{\Omega}\circ\varphi^{-1})'\in C(\varphi({\mathcal{U}}\cap\mathcal{U}_{\Omega}),
  L({C}^{m-1+\mu}(S),{C}^{m-1+\mu}(\Gamma)))$, where the second $\varphi({\mathcal{U}}\cap\mathcal{U}_{\Omega})$
  uses the topology of ${C}^{m+\mu}(S)$, it follows that $\varphi_{\Omega}\circ\varphi^{-1}$ is also continuously
  differentiable as a map from $\varphi({\mathcal{U}}\cap\mathcal{U}_{\Omega})\subseteq {C}^{m+\mu}(S)$ to
  ${C}^{m+\mu-1}(\Gamma)$, and $(\varphi_{\Omega}\circ\varphi^{-1})'\in C(\varphi({\mathcal{U}}\cap\mathcal{U}_{\Omega}),
  L({C}^{m-1+\mu}(S),{C}^{m-1+\mu}(\Gamma)))$. Since $F\circ\varphi_{\Omega}^{-1}$ is continuously differentiable in
  $O=\varphi_{\Omega}(\mathcal{O}'\cap\mathcal{U})$ in the topology of ${C}^{m-1+\mu}(\Gamma)$ and
  $(F\circ\varphi_{\Omega}^{-1})'\in C(O,L({C}^{m-1+\mu}(\Gamma),\mathbb{R}))$, where the second $O$ uses the topology of
  ${C}^{m+\mu}(\Gamma)$, from the above relation we see that $F\circ\varphi^{-1}$ is continuously differentiable in
  $\varphi(\mathcal{O}'\cap\mathcal{U})$ in the topology of ${C}^{m+\mu}(S)$, and $(F\circ\varphi^{-1})'\in
  C(\varphi(\mathcal{O}'\cap\mathcal{U}),L({C}^{m-1+\mu}(S),\mathbb{R}))\subseteq C(\varphi(\mathcal{O}'\cap\mathcal{U}),
  L({C}^{m+\mu}(S),\mathbb{R}))$, where $\varphi(\mathcal{O}'\cap\mathcal{U})$
  uses the topology of ${C}^{m+\mu}(S)$. Hence $F$ is fully strongly continuously differentiable at $\Omega$. This proves
  the assertion (1).

  (2)\ Assume that $f$ is differentiable at $\Omega$. Then there exists regular local chart $(\mathcal{U},\varphi)$ of
  ${\mathfrak{D}}^{m+\mu}({\mathbb{R}}^n)$ at $\Omega$, such that the function $t\mapsto\varphi(f(t))$ is differentiable
  at $t=0$ in the topology of $X={C}^{m+\mu}(S)$. We have
$$
  \varphi_{\Omega}(f(t))=(\varphi_{\Omega}\circ\varphi^{-1})(\varphi(f(t))).
$$
  As before, $\varphi_{\Omega}\circ\varphi^{-1}$ is continuously differentiable as a map from $\varphi({\mathcal{U}}\cap
  \mathcal{U}_{\Omega})\subseteq {C}^{m+\mu}(S)$ to ${C}^{m+\mu-1}(\Gamma)$. Hence, from the above relation we see
  that the function $t\mapsto\varphi_{\Omega}(f(t))$ is differentiable at $t=0$ in the topology of ${C}^{m-1+\mu}(\Gamma)$.
  Moreover, since $(\varphi\circ f)'(0)\in{C}^{m+\mu}(S)$, $(\varphi_{\Omega}\circ\varphi^{-1})'(\varphi(\Omega))\in
  L({C}^{m+\mu}(S),{C}^{m+\mu}(\Gamma))$ (by (4.15)), we have
$$
  (\varphi_{\Omega}\circ f)'(0)=(\varphi_{\Omega}\circ\varphi^{-1})'(\varphi(\Omega))(\varphi\circ f)'(0)\in{C}^{m+\mu}(\Gamma).
$$
  Hence $f$ is weakly differentiable at $\Omega$.

  (3)\ We denote by $A$ the set of the right-hand side of (4.17). By the assertion (2), it is clear that
  $T_{\Omega}({\mathfrak{D}}^{m+\mu}({\mathbb{R}}^n))\subseteq A$. We now prove the inverse inclusion. Let
  $u\in A$. Then there exists a curve $f:(-\varepsilon,\varepsilon)\to{\mathfrak{D}}^{m+\mu}({\mathbb{R}}^n)$
  $(\varepsilon>0)$ in ${\mathfrak{D}}^{m+\mu}({\mathbb{R}}^n)$ passing $\Omega$, i.e., $f(0)=\Omega$, $f(t)$
  is weakly differentiable at $t=0$, such that $f'(0)=u$. Let $\xi=(\varphi_{\Omega}\circ f)'(0)$ and
  $\eta=(\varphi\circ\varphi_{\Omega}^{-1})'(0)\xi$. Then $\xi\in{C}^{m+\mu}(\Gamma)$ and, by the second relation
  in (4.15), $\eta\in{C}^{m+\mu}(S)$. Moreover, $(\varphi_{\Omega}\circ\varphi^{-1})'(\varphi(\Omega))\eta=\xi$.
  Define a curve $g$ in ${\mathfrak{D}}^{m+\mu}({\mathbb{R}}^n)$ by letting $g(t)=\varphi^{-1}(t\eta)$ for
  $|t|<\varepsilon'$, where $\varepsilon'>0$ is sufficiently small. Clearly $g$ is differentiable at $t=0$, and
  $(\varphi\circ g)'(0)=\eta$. Hence, for any $F\in\hat{\mathscr{D}}^{1\!s}_{\Omega}$ we have
$$
\begin{array}{rl}
  g'(0)F=&(F\circ g)'(0)=(F\circ\varphi^{-1})'(\varphi(\Omega))(\varphi\circ g)'(0)
\\
  =&(F\circ\varphi_{\Omega}^{-1})'(0)(\varphi_{\Omega}\circ\varphi^{-1})'(\varphi(\Omega))\eta
  =(F\circ\varphi_{\Omega}^{-1})'(0)\xi
\\
  =&(F\circ\varphi_{\Omega}^{-1})'(0)(\varphi_{\Omega}\circ f)'(0)
  =(F\circ f)'(0)=f'(0)F.
\end{array}
$$
  Hence $g'(0)=u$. This proves $u\in T_{\Omega}({\mathfrak{D}}^{m+\mu}({\mathbb{R}}^n))$. Hence
  $T_{\Omega}({\mathfrak{D}}^{m+\mu}({\mathbb{R}}^n))=A$, and proves the assertion (3). The proof of Lemma 4.5
  is complete. $\quad\Box$
\medskip

  By Lemma 4.5, it follows that the following relation holds:
\begin{equation}
  T_{\Omega}({\mathfrak{D}}^{m+\mu}({\mathbb{R}}^n))\cong {C}^{m+\mu}(\Gamma).
\end{equation}
  Hence, the tangent space $T_{\Omega}({\mathfrak{D}}^{m+\mu}({\mathbb{R}}^n))$ of ${\mathfrak{D}}^{m+\mu}({\mathbb{R}}^n)$
  at a point $\Omega\in{\mathfrak{D}}^{m+\mu+1}({\mathbb{R}}^n)$ can be alternatively defined to be the Banach space
  ${C}^{m+\mu}(\Gamma)$: $T_{\Omega}({\mathfrak{D}}^{m+\mu}({\mathbb{R}}^n))={C}^{m+\mu}(\Gamma)$, but for this
  purpose we need to identify any $Q\in\mathcal{U}_{\Omega}$ with $\varphi_{\Omega}(Q)$ ($f'(0)$ is then identified with
  $(\varphi_{\Omega}\circ f)'(0)$). Note that in this definition the tangent vector $f'(0)$ of a curve $f:(-\varepsilon,\varepsilon)
  \to{\mathfrak{D}}^{m+\mu}({\mathbb{R}}^n)$ $(\varepsilon>0)$ has a simple physical explanation: Regarding $f$ as a flow of
  simple $C^{m+\mu}$-domains,  $f'(0)$ is the normal velocity of the boundary $\partial\Omega=\Gamma$. To see this let
  $f:(-\varepsilon,\varepsilon)\to{\mathfrak{D}}^{m+\mu}({\mathbb{R}}^n)$ be a curve in ${\mathfrak{D}}^{m+\mu}({\mathbb{R}}^n)$
  passing $\Omega$, differentiable at $t=0$. Let $f(t)=\Omega_t$, $\forall t\in (-\varepsilon,\varepsilon)$ (so that $\Omega_0=\Omega$).
  There exists a function $\rho:(-\varepsilon,\varepsilon)\times\Gamma\to{\mathbb{R}}$ such that
$$
  \partial\Omega_t=\{x+\rho(t,x)\bfn(x):x\in\Gamma\}, \quad \forall t\in (-\varepsilon,\varepsilon); \quad \rho(0,\cdot)=0.
$$
  Clearly, $(\varphi_{\Omega}\circ f)(t)=\rho(t,\cdot)$, $\forall t\in (-\varepsilon,\varepsilon)$. Hence $(\varphi_{\Omega}
  \circ f)'(0)=\partial_t\rho(0,\cdot)$, which is clearly the normal velocity of $\partial\Omega$. This is exactly the approach
  of \cite{PruSim} to define the tangent space $T_{\Omega}({\mathfrak{D}}^{m+\mu}({\mathbb{R}}^n))$.

  Note that a similar discussion can be made for $\dot{\mathfrak{D}}^{m+\mu}({\mathbb{R}}^n)$ and
  $T_{\Omega}(\dot{\mathfrak{D}}^{m+\mu}({\mathbb{R}}^n))$, which we omit here.

\section{Vector bundles over $\mathfrak{D}^{m+\mu}({\mathbb{R}}^n)$}
\setcounter{equation}{0}

\hskip 2em
  Next let us consider Banach vector bundles over the Banach manifold $\mathfrak{D}^{m+\mu}({\mathbb{R}}^n)$.
\medskip

  {\bf Definition 5.1}\ \ {\em Let $m,k$ be positive integers and $0\leqslant\mu,\nu\leqslant 1$. We denote by
  $\mathfrak{E}^{m+\mu,k+\nu}({\mathbb{R}}^n)$ the following Banach vector bundle over the Banach manifold
  $\mathfrak{D}^{m+\mu}({\mathbb{R}}^n)$:
$$
  \mathfrak{E}^{m+\mu,k+\nu}({\mathbb{R}}^n)=\{(\Omega,u):\;\Omega\in\mathfrak{D}^{m+\mu}({\mathbb{R}}^n),
  u\in C^{k+\nu}(\overline{\Omega})\}.
$$
  The fibre at the point $\Omega\in\mathfrak{D}^{m+\mu}({\mathbb{R}}^n)$ is the Banach space $C^{k+\nu}(\overline{\Omega})$
  which is isomorphic to $C^{k+\nu}(\bar{\mathbb{B}}(0,1))$, and the project is}
$$
  p(u,\Omega)=\Omega, \quad \forall(u,\Omega)\in\mathfrak{E}^{m+\mu,k+\nu}({\mathbb{R}}^n).
$$

  Surely, $\mathfrak{E}^{m+\mu,k+\nu}({\mathbb{R}}^n)$ is also a Banach manifold. Local chart of $\mathfrak{E}^{m+\mu,k+\nu}({\mathbb{R}}^n)$
  at a point $(\Omega,u)$ is as follows: Let $(\mathcal{U},\varphi)$ be a regular local chart of $\mathfrak{D}^{m+\mu}({\mathbb{R}}^n)$
  at the point $\Omega$ with base hypersurface $S$. Let $Q$ be the simple domain enclosed by $S$. Let $\mathcal{O}=
  \varphi(\mathcal{U})\in V^{m+\mu}(S)$ (cf. (3.1), (3.2)). For any $\rho\in\mathcal{O}$ let $\theta_{\rho}$ be as in $(3.1)'$.
  We use the Hanzawa transformation to extend $\theta_{\rho}:S\to\partial\Omega_{\rho}$ into a homeomorphism $\Theta_{\rho}:
  \overline{Q}\to\overline{\Omega}_{\rho}$ (cf. \cite{Cui5, EscSim2, EscSim3}). Note that $\Theta_{\rho}$ depends on $\rho$
  continuously. Recall that the pull-back $\Theta_{\rho}^*:C^{k+\nu}(\overline{\Omega}_{\rho})\to C^{k+\nu}(\overline{Q})$ induced
  by $\Theta_{\rho}$ is defined as follows:
$$
  \Theta_{\rho}^*(v)=v\circ\Theta_{\rho}, \quad \forall v\in C^{k+\nu}(\overline{\Omega}_{\rho}).
$$
  Now let
$$
  \widetilde{\mathcal{U}}=\{(R,v):\;R\in\mathcal{U},v\in C^{k+\nu}(\overline{R})\}
  =\{(\Omega_{\rho},v):\;\rho\in\mathcal{O},v\in C^{k+\nu}(\overline{\Omega}_{\rho})\},
$$
  and let $\widetilde{\varphi}:\widetilde{\mathcal{U}}\to C^{m+\mu}(\overline{Q})\times C^{k+\nu}(\overline{Q})$ be the following map:
$$
  \widetilde{\varphi}(\Omega_{\rho},v)=(\rho,\Theta_{\rho}^*(v)), \quad
  \forall\rho\in\mathcal{O},\;\;\forall v\in C^{k+\nu}(\overline{\Omega}_{\rho}).
$$
  Then $(\widetilde{\mathcal{U}},\widetilde{\varphi})$ is a local chart of $\mathfrak{E}^{m+\mu,k+\nu}({\mathbb{R}}^n)$ at the point
  $(\Omega,u)$.
\medskip

  {\bf Lemma 5.2}\ \ {\em Let $(\mathcal{U}_1,\varphi_1)$, $(\mathcal{U}_2,\varphi_2)$ be two local charts of
  $\mathfrak{D}^{m+\mu}({\mathbb{R}}^n)$ such that $\mathcal{U}_1\cap\mathcal{U}_2\neq\varnothing$. Let
  $(\widetilde{\mathcal{U}}_1,\widetilde{\varphi}_1)$, $(\widetilde{\mathcal{U}}_2,\widetilde{\varphi}_2)$ be the
  corresponding local charts of $\mathfrak{E}^{m+\mu,k+\nu}({\mathbb{R}}^n)$ as defined above for a given point
  $\Omega\in\mathcal{U}_1\cap\mathcal{U}_2$. If $\Omega\in\mathcal{U}_1\cap\mathcal{U}_2\cap
  \mathfrak{D}^{m+\mu+1}({\mathbb{R}}^n)$ and $z\in C^{k+\nu+1}(\overline{Q}_1)$, then $\widetilde{\varphi}_2\circ
  \widetilde{\varphi}_1^{-1}$ is differentiable at $(\rho,z)$, where $\rho=\varphi_1(\Omega)$. Moreover,
  $[(\rho,z)\mapsto(\widetilde{\varphi}_2\circ\widetilde{\varphi}_1^{-1})'(\rho,z)]\in C((\mathcal{O}_1\cap C^{m+\mu+1}(S_1))
  \times C^{k+\nu+1}(\overline{Q}_1),L(C^{m+\mu}(S_1)\times C^{k+\nu}(\overline{Q}_1),
  C^{m+\mu}(S_2)\times C^{k+\nu}(\overline{Q}_2)))$, where $Q_1$, $Q_2$ are the domains enclosed by the base hypersurfaces
  $S_1$, $S_2$ of the local charts $(\mathcal{U}_1,\varphi_1)$, $(\mathcal{U}_2,\varphi_2)$, respectively.}
\medskip

  {\em Proof}.\ \ Given $R\in\mathcal{U}_1\cap\mathcal{U}_2$, we denote $\rho_i=\varphi_i(R)\in C^{m+\mu}(S_i)$, $i=1,2$,
  and we write $R=\Omega^i_{\rho_i}$, $i=1,2$. Also we use the notation $\Theta^i_{\rho_i}$ to denote the Hanzawa
  transformation related to the map (see $(3.1)'$) $\theta^i_{\rho}:S_i\to\partial R=\partial\Omega^i_{\rho_i}$, $i=1,2$,
  so that $\Theta^i_{\rho_i}:\overline{Q}_i\to\overline{R}=\overline{\Omega}^i_{\rho_i}$, $i=1,2$. It follows that
$$
  \widetilde{\varphi}_i(R,v)=(\rho_i,\Theta_{\rho_i}^{i*}(v)), \quad \forall R\in\mathcal{U}_1\cap\mathcal{U}_2,\;\;
  \forall v\in C^{k+\nu}(\overline{R}),\;\; i=1,2.
$$
  Hence,
$$
  (\widetilde{\varphi}_2\circ\widetilde{\varphi}_1^{-1})(\rho,z)
  =((\varphi_2\circ\varphi_1^{-1})(\rho),(\Theta_{\rho_2}^{2*}\circ\Theta_{\rho_1}^{1*-1})(z)),
  \quad \forall \rho\in\mathcal{O}_1=\varphi_1(\mathcal{U}_1\cap\mathcal{U}_2),\;\;
  \forall z\in C^{k+\nu}(\overline{Q}_1),
$$
  where $\rho_1=\rho$, $\rho_2=(\varphi_2\circ\varphi_1^{-1})(\rho)$. Note that
$$
  (\Theta_{\rho_2}^{2*}\circ\Theta_{\rho_1}^{1*-1})(z)=z\circ(\Theta_{\rho_1}^1)^{-1}\circ\Theta_{\rho_2}^2.
$$
  From these expressions, we easily see that if $\Omega\in\mathcal{U}_1\cap\mathcal{U}_2\cap
  \mathfrak{D}^{m+\mu+1}({\mathbb{R}}^n)$ and $z\in C^{k+\nu+1}(\overline{Q}_1)$, then $\widetilde{\varphi}_2
  \circ\widetilde{\varphi}_1^{-1}$ is differentiable at $(\rho,z)$, and $[(\rho,z)\mapsto(\widetilde{\varphi}_2
  \circ\widetilde{\varphi}_1^{-1})'(\rho,z)]\in C((\mathcal{O}_1\cap C^{m+\mu+1}(S_1))\times
  C^{k+\nu+1}(\overline{Q}_1),L(C^{m+\mu}(S_1)\times C^{k+\nu}(\overline{Q}_1),C^{m+\mu}(S_2)\times
  C^{k+\nu}(\overline{Q}_2)))$. $\quad\Box$
\medskip

  In view of the above lemma, we can call points of $\mathfrak{E}^{m+\mu+1,k+\nu+1}({\mathbb{R}}^n)$ as {\em differentiable
  points} of $\mathfrak{E}^{m+\mu,k+\nu}({\mathbb{R}}^n)$, and similarly as in Definition define a tangent space of
  $\mathfrak{E}^{m+\mu,k+\nu}({\mathbb{R}}^n)$ at every point $(\Omega,u)\in\mathfrak{E}^{m+\mu+1,k+\nu+1}({\mathbb{R}}^n)$,
  which we denote as $T_{(\Omega,u)}(\mathfrak{E}^{m+\mu,k+\nu}({\mathbb{R}}^n))$. Note that
$$
  T_{(\Omega,u)}(\mathfrak{E}^{m+\mu,k+\nu}({\mathbb{R}}^n))\cong
  C^{m+\mu}(\mathbb{S}^{n-1})\times C^{k+\nu}(\overline{\mathbb{B}^n}).
$$

  Note that a similar discussion can be made for $\dot{\mathfrak{D}}^{m+\mu}({\mathbb{R}}^n)$, which we omit here. The
  corresponding vector bundle is denoted as $\dot{\mathfrak{E}}^{m+\mu,k+\nu}({\mathbb{R}}^n)$.

\section{Lie group actions to $\mathfrak{D}^{m+\mu}({\mathbb{R}}^n)$}
\setcounter{equation}{0}

\hskip 2em
  It is clear that $\mathfrak{D}^{m+\mu}({\mathbb{R}}^n)$ is invariant under dilations, translations and rotations
  in ${\mathbb{R}}^n$. In what follows we study smoothness of these Lie group actions to $\mathfrak{D}^{m+\mu}({\mathbb{R}}^n)$.

  We first consider the action of the translation group. Let $G_{tl}={\mathbb{R}}^n$ be the additive group of $n$-vectors.
  Given $z\in{\mathbb{R}}^n$ and $\Omega\in\mathfrak{D}^{m+\mu}({\mathbb{R}}^n)$, let
$$
  p(z,\Omega)=\Omega+z=\{x+z:\,x\in\Omega\}.
$$
  It is clear that $p(z,\Omega)\in\mathfrak{D}^{m+\mu}({\mathbb{R}}^n)$, $\forall\Omega\in\mathfrak{D}^{m+\mu}({\mathbb{R}}^n)$,
  $\forall z\in{\mathbb{R}}^n$. Moreover, it is also clear that
$$
  p(0,\Omega)=\Omega, \quad p(z_1+z_2,\Omega)=p(z_1,p(z_2,\Omega)), \quad
  \forall\Omega\in\mathfrak{D}^{m+\mu}({\mathbb{R}}^n), \;\; \forall z_1,z_2\in
  {\mathbb{R}}^n.
$$
  Hence $(G_{tl},p)$ is a Lie group action on $\mathfrak{D}^{m+\mu}({\mathbb{R}}^n)$.
\medskip

  {\bf Lemma 6.1}\ \ {\em The Lie group action $(G_{tl},p)$ on $\mathfrak{D}^{m+\mu}({\mathbb{R}}^n)$ satisfies the following
  properties: For any nonnegative integers $k$ and $l$,
\begin{equation}
  [z\mapsto p(z,\cdot)]\in C^k(G_{tl},\mathfrak{C}^{\,l}(\mathfrak{D}^{m+k+l+\mu}({\mathbb{R}}^n),\mathfrak{D}^{m+\mu}({\mathbb{R}}^n))).
\end{equation}
  In particular, $p\in C(G_{tl}\times\mathfrak{D}^{m+\mu}({\mathbb{R}}^n),\mathfrak{D}^{m+\mu}({\mathbb{R}}^n))$.}
\medskip

  {\em Proof}.\ \ We first compute the representation of the mapping $p:G_{tl}\times\mathfrak{D}^{m+\mu}({\mathbb{R}}^n)\to
  \mathfrak{D}^{m+\mu}({\mathbb{R}}^n)$ in local charts of $G_{tl}\times\mathfrak{D}^{m+\mu}({\mathbb{R}}^n)$ and
  $\mathfrak{D}^{m+\mu}({\mathbb{R}}^n)$. Let $z_0\in G_{tl}$ and $\Omega\in\mathfrak{D}^{m+\mu}({\mathbb{R}}^n)$ be given.
  Choose a closed $C^{\infty}$ hypersurface $S\subseteq{\mathbb{R}}^n$ sufficiently closed to $\partial\Omega$ such that
  it satisfies the conditions $(a)$--$(c)$ in the previous subsection. Let $(\mathcal{U},\varphi)$ be the local chart of
  $\mathfrak{D}^{m+\mu}({\mathbb{R}}^n)$ at the point $\Omega$ as defined by (3.6) and (3.7). Let
$$
  \hat{\Omega}=\Omega+z_0, \quad  \hat{S}=S+z_0=\{x+z_0:x\in S\},
$$
  and $(\hat{\mathcal{U}},\hat{\varphi})$ be the local chart of $\mathfrak{D}^{m+\mu}({\mathbb{R}}^n)$ at the point $\hat{\Omega}$
  as defined by (3.6) and (3.7), with $S$ there replaced by $\hat{S}$. Take $\varepsilon>0$ sufficiently small such that,
  by slightly shrinking the neighborhood $\mathcal{U}$ of $\Omega$ when necessary, we have $Q+z\in\hat{\mathcal{U}}$ for
  all $Q\in\mathcal{U}$ and $z\in B(z_0,\varepsilon)$. Here $B(z_0,\varepsilon)$ denotes the open sphere in $G_{tl}=
  {\mathbb{R}}^n$ with center $z_0$ and radius $\varepsilon$. Now let $Q\in\mathcal{U}$ and $z\in B(z_0,\varepsilon)$.
  There exists a unique $\rho\in\mathcal{O}$, where $\mathcal{O}$ is a neighborhood of the origin in $C^{m+\mu}(S)$, such
  that
\begin{equation}
  \partial Q=\{x+\rho(x)\bfn(x):x\in S\},
\end{equation}
  where $\bfn$ is the outward unit normal field of $S$. By definition we have $\varphi(Q)=\rho$. Similarly, there exists
  a unique $\hat{\rho}\in\hat{\mathcal{O}}$, where $\hat{\mathcal{O}}$ is a neighborhood of the origin in $C^{m+\mu}(\hat{S})$,
  such that
\begin{equation}
  \partial(Q+z)=\{x+\hat{\rho}(x)\hat{\bfn}(x):x\in\hat{S}\},
\end{equation}
  where $\hat{\bfn}$ is the outward unit normal field of $\hat{S}$. By definition we have $\hat{\varphi}(Q+z)=\hat{\rho}$.
  Hence, letting
$$
  \tilde{p}(z,\rho)=(\hat{\varphi}\circ p)(z,\varphi^{-1}(\rho)), \quad \forall\rho\in\mathcal{O},\;\; \forall z\in B(z_0,\varepsilon),
$$
  we have
$$
  \tilde{p}(z,\rho)=\hat{\rho}, \quad \forall\rho\in\mathcal{O},\;\; \forall z\in B(z_0,\varepsilon).
$$
  Next, from (6.2) we have
$$
  \partial(Q+z)=\{x+z+\rho(x)\bfn(x):x\in S\}=\{y-z_0+z+\rho(y-z_0)\bfn(y-z_0):y\in\hat{S}\}.
$$
  Comparing this expression with (4.3) we get
\begin{equation}
  x+\hat{\rho}(x)\hat{\bfn}(x)=y+z-z_0+\rho(y-z_0)\bfn(y-z_0), \quad x,y\in\hat{S},
\end{equation}
  which, similarly as before, means that for any $x\in\hat{S}$ there exists a unique $y\in\hat{S}$ such that the above equality
  holds, and vice versa. It follows that there exists a function $y=\mu(\rho,z,x)$ uniquely determined by $(\rho,z)$, mapping
  $x\in\hat{S}$ to $y\in\hat{S}$, such that
\begin{equation}
  \hat{\rho}(x)=\langle\mu(\rho,z,x)-x+z-z_0,\hat{\bfn}(x)\rangle+\rho(\mu(\rho,z,x)-z_0)\langle\bfn(\mu(\rho,z,x)-z_0),\hat{\bfn}(x)\rangle,
  \quad x\in\hat{S}.
\end{equation}
  Similarly as in the proof of Theorem 3.4 (cf. (3.14)), we have
$$
  [(\rho,z)\mapsto\mu(\rho,z,\cdot)]\in C^k(\mathcal{O}\times B(z_0,\varepsilon),C^{m-k+\mu}(\hat{S},\hat{S})), \quad k=0,1,\cdots,m.
$$
  Using this fact, the expression (6.5) and a similar argument as in the proof of Lemma 3.5 we get
$$
\begin{array}{rl}
  \{z\mapsto[\rho\mapsto\partial_{\rho}^j\tilde{p}(z,\rho)]\}\in &C^k(B(z_0,\varepsilon),C^{\,l-j}(\mathcal{O}\cap C^{m+k+l+\mu}(S),
  L^j(C^{m+\mu}(S),C^{m+\mu}(\hat{S})))),\\
  & \qquad\qquad\qquad  k,l=0,1,\cdots,\;\; j=0,1,\cdots,l.
\end{array}
$$
  This proves (6.1). $\quad\Box$
\medskip

  Next we consider the action of the dilation group. Let $G_{\!dl}={\mathbb{R}}_+=(0,\infty)$ be the multiplicative group of
  all positive numbers. Given $\lambda\in G_{\!dl}$ and $\Omega\in\mathfrak{D}^{m+\mu}({\mathbb{R}}^n)$, let
$$
  q(\lambda,\Omega)=\lambda\Omega=\{\lambda x:\,x\in\Omega\}.
$$
  Clearly $q(\lambda,\Omega)\in\mathfrak{D}^{m+\mu}({\mathbb{R}}^n)$, $\forall\Omega\in\mathfrak{D}^{m+\mu}({\mathbb{R}}^n)$,
  $\forall\lambda\in G_{\!dl}$. Moreover, it is also clear that
$$
  q(1,\Omega)=\Omega, \quad q(\lambda_1\lambda_2,\Omega)=q(\lambda_1,q(\lambda_2,\Omega)), \quad
  \forall\Omega\in\mathfrak{D}^{m+\mu}({\mathbb{R}}^n), \;\; \forall\lambda_1,\lambda_2\in G_{\!dl}.
$$
  Hence $(G_{\!dl},p)$ is a Lie group action on $\mathfrak{D}^{m+\mu}({\mathbb{R}}^n)$. Similar to Lemma 6.1 we have
\medskip

  {\bf Lemma 6.2}\ \ {\em The Lie group action $(G_{\!dl},q)$ in $\mathfrak{D}^{m+\mu}({\mathbb{R}}^n)$ satisfies the
  following properties: For any nonnegative integers $k$ and $l$,
$$
  [\lambda\mapsto q(\lambda,\cdot)]\in C^k(G_{\!dl},\mathfrak{C}^{\,l}(\mathfrak{D}^{m+k+l+\mu}({\mathbb{R}}^n),
  \mathfrak{D}^{m+\mu}({\mathbb{R}}^n))).
$$
  In particular, $q\in C(G_{\!dl}\times\mathfrak{D}^{m+\mu}({\mathbb{R}}^n),\mathfrak{D}^{m+\mu}({\mathbb{R}}^n))$.}
\medskip

  The proof is similar to that of Lemma 6.1; we omit it. $\quad\Box$
\medskip

  The group actions $(G_{tl},p)$ and $(G_{\!dl},q)$ to $\mathfrak{D}^{m+\mu}({\mathbb{R}}^n)$ are not mutually commutative.
  However, it is clear that the following relation holds:
\begin{equation}
  q(\lambda,p(z,\Omega))=p(\lambda z, q(\lambda,\Omega)), \quad \forall\Omega\in\mathfrak{D}^{m+\mu}({\mathbb{R}}^n),\;\;
  \forall z\in G_{tl},\;\; \forall\lambda\in G_{\!dl}.
\end{equation}
  Hence the actions $p$ of $G_{tl}$ and $q$ of $G_{\!dl}$ to $\mathfrak{D}^{m+\mu}({\mathbb{R}}^n)$ are quasi-commutative. Besides, denoting
$$
  g(z,\lambda,\Omega)=p(z, q(\lambda,\Omega)), \quad \forall\Omega\in\mathfrak{D}^{m+\mu}({\mathbb{R}}^n),\;\;
  \forall z\in G_{tl},\;\; \forall\lambda\in G_{\!dl},
$$
  we easily see that the following relation holds:
\begin{equation}
  {\rm rank}\,\partial_{(z,\lambda)}g(z,\lambda,\Omega)=n+1, \quad \forall\Omega\in\mathfrak{D}^{m+1+\mu}({\mathbb{R}}^n),\;\;
  \forall z\in G_{tl},\;\; \forall\lambda\in G_{\!dl}.
\end{equation}

  Finally, let $O(n)$ be the Lie group of all $n\times n$ orthogonal matrices. For any $A\in O(n)$ and $\Omega\in
  \mathfrak{D}^{m+\mu}({\mathbb{R}}^n)$, let
$$
  r(A,\Omega)=\mathbf{A}(\Omega)=\{\mathbf{A}(x):\,x\in\Omega\},
$$
  where $\mathbf{A}$ denotes the orthogonal transformation in ${\mathbb{R}}^n$ induced by $A$, i.e, $\mathbf{A}(x)=(Ax^T)^T$ by
  regarding vectors in ${\mathbb{R}}^n$ as $1\times n$ matrix. Clearly $r(A,\Omega)\in\mathfrak{D}^{m+\mu}({\mathbb{R}}^n)$,
  $\forall\Omega\in\mathfrak{D}^{m+\mu}({\mathbb{R}}^n)$, $\forall A\in O(n)$. Moreover, it is also
  clear that
$$
  r(I,\Omega)=\Omega, \quad r(A_1A_2,\Omega)=r(A_1,r(A_2,\Omega)), \quad \forall\Omega\in\mathfrak{D}^{m+\mu}({\mathbb{R}}^n), \;\;
  \forall A_1,A_2\in O(n).
$$
  Hence $(O(n),r)$ is a Lie group action on $\mathfrak{D}^{m+\mu}({\mathbb{R}}^n)$. Similar to Lemmas 6.1 and 6.2 we have
\medskip

  {\bf Lemma 6.3}\ \ {\em Let $m,n\geqslant 2$ be integers and $0\leqslant\mu\leqslant 1$. The Lie group action $(O(n),q)$ in
  $\mathfrak{D}^{m+\mu}({\mathbb{R}}^n)$ satisfies the following properties: For any nonnegative integers $k$ and $l$,
$$
  [A\mapsto r(A,\cdot)]\in C^k(O(n),\mathfrak{C}^{\,l}(\mathfrak{D}^{m+k+l+\mu}({\mathbb{R}}^n),\mathfrak{D}^{m+\mu}({\mathbb{R}}^n))).
$$
  In particular, $r\in C(O(n)\times\mathfrak{D}^{m+\mu}({\mathbb{R}}^n),\mathfrak{D}^{m+\mu}({\mathbb{R}}^n))$.}
\medskip

  The proof is similar to that of Lemma 6.1; we omit it. $\quad\Box$

  As we shall not use the group action $(O(n),r)$ later on, we do not make further discussion to it here.

  Note that all discussions made in this section apply to $\dot{\mathfrak{D}}^{m+\mu}({\mathbb{R}}^n)$ when all $C^{m+\mu}$-spaces
  are replaced with corresponding $\dot{C}^{m+\mu}$-spaces. Here we do not repeat the details. Later on we shall use such results
  without further discussion.

  Note that a similar discussion can be made for $\dot{\mathfrak{D}}^{m+\mu}({\mathbb{R}}^n)$, which we omit here.

\section{Applications to free boundary problems}
\setcounter{equation}{0}

\hskip 2em
  In this section we give some examples to show how to reduce free boundary problems into differential equations
  in the Banach manifold $\dot{\mathfrak{D}}^{m+\mu}(\mathbb{R}^n)$.

  Let us first consider the one phase Hele-Shaw problem in the whole space ${\mathbb{R}}^n$:
\begin{equation}
\left\{
\begin{array}{ll}
   \Delta u(x,t)=0, &\qquad x\in\Omega(t), \;\; t>0,\\
    u(x,t)=-\kappa(x,t),   &\qquad x\in\partial\Omega(t), \;\; t>0,\\
    V_n(x,t)=\partial_{\bfn}u(x,t), &\qquad x\in\partial\Omega(t), \;\; t>0,\\
   \Omega(0)=\Omega_0, &
\end{array}
\right.
\end{equation}
  where for each $t>0$, $\Omega(t)$ is an unknown domain in ${\mathbb{R}}^n$, $\Delta$ is the Laplacian in $n$ variables,
  $u=u(x,t)$ is an unknown function defined for $x\in\overline{\Omega(t)}$ and $t\geqslant 0$, $\kappa(\cdot,t)$ is the
  mean curvature of the boundary $\partial\Omega(t)$ of $\Omega(t)$, $V_n$ is the normal velocity of the free boundary
  $\partial\Omega(t)$, $\bfn$ denotes the outward-pointing normal field of $\partial\Omega(t)$, and $\Omega_0$ is a given
  initial domain. Here we take the convention that for a convex domain the mean curvature of its boundary takes nonnegative
  values. For simplicity we only consider the case that $\Omega_0$ is a simple domain sufficiently close to a sphere.

  The above problem has been intensively studied during the past fifty years, cf. \cite{Chen}, \cite{ChenHY}, \cite{EscSim2},
  \cite{EscSim3}, \cite{EscSim4}, \cite{EscSim5} and references therein. The purpose here is to show that this problem can be
  rewritten as a differential equation in the Banach manifold $\mathfrak{M}:=\dot{\mathfrak{D}}^{m+\mu}({\mathbb{R}}^n)$ for
  any integer $m\geqslant 1$ and $0<\mu<1$. Let $\mathfrak{M}_0:=\dot{\mathfrak{D}}^{m+3+\mu}({\mathbb{R}}^n)$. We know that
  $\mathfrak{M}$ and $\mathfrak{M}_0$ are Banach manifolds built on the Banach spaces $\dot{C}^{m+\mu}({\mathbb{S}}^{n-1})$ and
  $\dot{C}^{m+3+\mu}({\mathbb{S}}^{n-1})$, respectively, and $\dot{C}^{m+3+\mu}({\mathbb{S}}^{n-1})$ is dense in
  $\dot{C}^{m+\mu}({\mathbb{S}}^{n-1})$.
  For any $Q\in\mathfrak{M}_0\subseteq\mathfrak{M}$ we use the standard local chart of $\mathfrak{M}$ at $Q$ to identify the
  tangent space $\mathcal{T}_{Q}(\mathfrak{M})$ with the Banach space $C^{m+\mu}(\partial Q)$. It follows that if
  $I\subseteq{\mathbb{R}}$ is an open interval and $\Omega:I\to\mathfrak{M}_0$ is a $C^1$ curve, then
$$
   \Omega'(t)=V_n(\cdot,t) \quad \mbox{for}\;\; t\in I.
$$
  Given $\Omega\in\mathfrak{M}_0$ we denote by $u_{\Omega}$ the unique solution of the following Dirichlet problem:
$$
\left\{
\begin{array}{ll}
   \Delta u_{\Omega}=0 &\quad\;\; \mbox{in} \;\; \Omega,\\
   u_{\Omega}=-\kappa_{\partial\Omega}   &\quad \mbox{on} \;\; \partial\Omega,
\end{array}
\right.
$$
  where $\kappa_{\partial\Omega}$ denotes the mean curvature of $\partial\Omega$. It is known that $u_{\Omega}\in
  \dot{C}^{m+\mu+1}(\overline{\Omega})$. Now introduce a vector field $\mathscr{F}$ in $\mathfrak{M}$ with domain $\mathfrak{M}_0$
  as follows: For any $\Omega\in\mathfrak{M}_0$ we define
$$
   \mathscr{F}(\Omega)=\partial_{\bfn}u_{\Omega}|_{\partial\Omega}\in\dot{C}^{m+\mu}(\partial\Omega)
   =\mathcal{T}_{\Omega}(\mathfrak{M}).
$$
  It follows that the problem (5.1) reduces into the following initial value problem of a differential equation in the Banach
  manifold $\mathfrak{M}$:
\begin{equation}
\left\{
\begin{array}{ll}
   \Omega'(t)=\mathscr{F}(\Omega(t)), &\quad  t>0,\\
    \Omega(0)=\Omega_0. &
\end{array}
\right.
\end{equation}
  Moreover, the vector field $\mathscr{F}$ is invariant under the translation group action $(G_{tl},p)$ and the rotation group
  $(O(n),r)$, i.e.,
\begin{equation}
  \mathscr{F}(p(a,\Omega))=\partial_{\Omega}p(a,\Omega)\mathscr{F}(\Omega), \quad \forall a\in G_{tl}, \;\; \forall\Omega\in\mathfrak{M}_0,
\end{equation}
\begin{equation}
  \mathscr{F}(r(A,\Omega))=\partial_{\Omega}r(A,\Omega)\mathscr{F}(\Omega), \quad \forall A\in O(n), \;\; \forall\Omega\in\mathfrak{M}_0,
\end{equation}
  and quasi-invariant under the dilation group action $(G_{dl},q)$ with quasi-invariant factor $\theta(\lambda)=\lambda^{-3}$,
  $\lambda>0$. i.e.,
\begin{equation}
  \mathscr{F}(q(\lambda,\Omega)))=\lambda^{-3}\partial_{\Omega}q(\lambda,\Omega)\mathscr{F}(\Omega), \quad \forall a\in G_{dl}, \;\;
  \forall\Omega\in\mathfrak{M}_0.
\end{equation}
  These assertions follow from the fact that the problem (7.1) is invariant under translation and rotation of coordinate,
  and quasi-invariant under dilation. Using these facts and applying some abstract results for parabolic differential
  equations in Banach manifolds \cite{Cui4}, we can re-establish the results obtained in the literatures mentioned above in a concise
  manner. We omit it here.

  Similar approach can be used to study the following free boundary problem modeling the growth of tumors:
\begin{equation}
\left\{
\begin{array}{ll}
\Delta\sigma=f(\sigma),&\quad x\in\Omega(t),\,\,t>0,\\
-\Delta p=g(\sigma),&\quad  x\in\Omega(t),\,\,t>0,\\
-\partial_{\mathbf{n}}\sigma=h(\sigma),&\quad  x\in\partial\Omega(t),\,\,t>0,\\
p=\gamma\kappa,&\quad  x\in\partial\Omega(t),\,\,t>0,\\
V_{\mathbf{n}}=-\partial_{\mathbf{n}}p,&\quad  x\in\partial\Omega(t),\,\,t>0,\\
\Omega(0)=\Omega_{0}.&
\end{array}\right.
\end{equation}
  Here $\Omega(t)$ is the domain in $\mathbb{R}^{n}$ occupied by the tumor at time $t$, $\sigma=\sigma(x,t)$ and $p=p(x,t)$
  are the nutrient concentration in the tumor region and the pressure between tumor cells, respectively,
  $\partial_{\mathbf{n}}$ represents the derivative in the direction of the outward normal $\mathbf{n}$ of the tumor
  surface $\partial\Omega(t)$, $\kappa$ is the mean curvature of the tumor surface $\partial\Omega(t)$ whose sign is designated
  by the convention that for the sphere it is positive, $\gamma$ is another positive constant reflecting the surface tension
  of the tumor surface and is usually referred to as {\em surface tension coefficient}, $V_{\mathrm{n}}$ is the normal velocity
  of the tumor surface movement, $f$, $g$ and $h$ are given functions with $f(\sigma)$ being the (normalized) consumption
  rate of nutrient by tumor cells when its concentration is at level $\sigma$, $g(\sigma)$ the (normalized) proliferation rate of
  tumor cells when the nutrient concentration is at level $\sigma$, and $h$ is a function measuring the combined effect of
  the strength of nutrient supply in its host tissue and the ability that the tumor receives nutrient from
  its surface, and $\Omega_{0}$ is the domain that the tumor initially occupies. Naturally, from physical viewpoint we have
  $n=3$; but for mathematical interest we consider the general case $n\geqslant 2$.

  Given $\Omega\in\mathfrak{M}_0:=\dot{\mathfrak{D}}^{m+3+\mu}({\mathbb{R}}^n)$, let us consider the following elliptic boundary
  value problem:
\begin{equation}
\left\{
\begin{array}{ll}
   \Delta\sigma=f(\sigma),&\quad \mbox{in}\;\;\Omega,\\
  -\Delta p=g(\sigma),&\quad \mbox{in}\;\;\Omega,\\
  -\partial_{\mathbf{n}}\sigma=h(\sigma),&\quad \mbox{on}\;\;\partial\Omega,\\
  p=\gamma\kappa_{\partial\Omega},&\quad \mbox{on}\;\;\partial\Omega,
\end{array}
\right.
\end{equation}
  where $\kappa_{\partial\Omega}$ denotes the mean curvature of $\partial\Omega$. Since $\partial\Omega$ is of $\dot{C}^{m+3+\mu}$-class
  and $f,h\in C^{\infty}[0,\infty)$, by applying the standard theory for elliptic boundary value problems we see that equation $(7.7)_1$
  subject to the boundary value condition $(7.7)_3$ has a unique solution $\sigma\in\dot{C}^{m+3+\mu}(\overline{\Omega})$ (existence
  follows from the fact that $\sigma=0$ and $\sigma=\bar{\sigma}$ are a pair of lower and upper solutions, and uniqueness follows
  from the condition $f',h'>0$). This implies that $g(\sigma)\in\dot{C}^{m+3+\mu}(\overline{\Omega})$.
  Since $\kappa_{\partial\Omega}\in\dot{C}^{m+1+\mu}(\partial\Omega)$, again by applying the standard theory for elliptic boundary
  value problems we see that equation $(7.7)_2$ subject to the boundary value condition $(7.7)_4$ has a unique solution
  $p\in\dot{C}^{m+1+\mu}(\overline{\Omega})$. Now we introduce a vector field $\mathscr{F}$ in $\mathfrak{M}$ with domain
  $\mathfrak{M}_0$ as follows: For any $\Omega\in\mathfrak{M}_0$ we define
$$
   \mathscr{F}(\Omega)=-\partial_{\bfn}p_{\Omega}|_{\partial\Omega}\in\dot{C}^{m+\mu}(\partial\Omega)
   =\mathcal{T}_{\Omega}(\mathfrak{M}).
$$
  It follows that the problem (7.6) reduces into the initial value problem (7.2) of a differential equation in the Banach
  manifold $\mathfrak{M}$. Again, since the problem (7.6) is invariant under translation and rotation of coordinate, it follows
  that the vector field $\mathscr{F}$ defined here also satisfies the relations (7.2) and (7.3). Using these facts and
  applying some abstract results for parabolic differential equations in Banach manifolds, we can prove that the problem
  (7.2) reduced from (7.6) is locally well-posed in $\dot{\mathfrak{D}}^{m+1+\mu}({\mathbb{R}}^n)$, and the the solution
  converges to a stationary solution when it initially is close to a stationary solution, cf. \cite{ZhengCui} for details.

  Finally let us consider the following free boundary problem modeling the motion of liquid drops (cf. \cite{GunPro, Sol1, Sol2}):
\begin{equation}
\left\{
\begin{array}{cl}
   \Delta\bfv=(\bfv\cdot\nabla)\bfv+\nabla\pi, &\quad\;\; x\in\Omega(t), \;\; t>0,\\
   \nabla\cdot\bfv=0,   &\quad \;\; x\in\Omega(t), \;\; t>0,\\
   T(\bfv,\pi)\bfn=-\gamma\kappa\bfn,   &\quad \;\; x\in\partial\Omega(t), \;\; t>0,\\
   V_n=\bfv\cdot\bfn, &\quad \;\; x\in\partial\Omega(t), \;\; t>0.
\end{array}
\right.
\end{equation}
  Here $\Omega(t)$ is the domain in ${\mathbb{R}}^3$ occupied by the liquid drop at time $t$, $\bfv$ and $\pi$ are velocity
  and pressure fields in the liquid drop, respectively, $\gamma>0$ is the surface tension coefficient constant, $\kappa$ is
  the mean curvature field on $\partial\Omega(t)$ whose sign is designated to be positive for spheres, $\bfn$ is the unit
  outward normal field on $\partial\Omega(t)$, $V_n$ is the normal velocity of the free boundary $\partial\Omega(t)$,
  $T(\bfv,\pi)=\nu S(\bfv)-\pi I$ ($I=$ the third-order identity matrix) is the stress tensor, and $S(\bfv)=\nabla\otimes\bfv
  +(\nabla\otimes\bfv)^T$ is the doubled strain tensor.

  The above problem and its more generalized forms have been extensively studied during the past fifty years, cf. the
  references therein. Here we show that the above problem can also be reduced into a differential equation in the Banach
  manifold $\dot{\mathfrak{D}}^{m+\mu}(\mathbb{R}^n)$. Indeed, given $\Omega\in\mathfrak{M}_0:=
  \dot{\mathfrak{D}}^{m+1+\mu}({\mathbb{R}}^n)$, let $(\bfv,\pi)$ be the solution of the following elliptic boundary value problem:
$$
\left\{
\begin{array}{cl}
   \Delta\bfv=(\bfv\cdot\nabla)\bfv+\nabla\pi, &\quad\;\; x\in\Omega,\\
   \nabla\cdot\bfv=0,   &\quad \;\; x\in\Omega,\\
   T(\bfv,\pi)\bfn=-\gamma\kappa\bfn,   &\quad \;\; x\in\partial\Omega.
\end{array}
\right.
$$
  It is well-known that this problem has a unique such solution $(\bfv,\pi)\in [\dot{C}^{m+\mu}(\overline{\Omega})]^3\times
  \dot{C}^{m-1+\mu}(\overline{\Omega})$. Now we introduce a vector field $\mathscr{F}$ in $\mathfrak{M}$ with domain
  $\mathfrak{M}_0$ as follows: For any $\Omega\in\mathfrak{M}_0$ we define
$$
   \mathscr{F}(\Omega)=\bfv\cdot\bfn|_{\partial\Omega}\in\dot{C}^{m+\mu}(\partial\Omega)
   =\mathcal{T}_{\Omega}(\mathfrak{M}).
$$
  It follows that the problem (7.8) reduces into the initial value problem (7.2) of a differential equation in the Banach
  manifold $\mathfrak{M}$.

  It is clear that (7.8) is invariant or quasi-invariant under translation and rotation of coordinate, and quasi-invariant
  under scaling of coordinate. Hence, the vector field $\mathscr{F}$ satisfies the relations (7.3)--(7.5). Using these facts
  and applying some abstract results for parabolic differential equations in Banach manifolds, we can prove that the problem
  (7.2) reduced from (7.8) is locally well-posed in $\dot{\mathfrak{D}}^{m+1+\mu}({\mathbb{R}}^n)$, and the the solution
  converges to a stationary solution when it initially is close to a stationary solution. These results have already been
  proved in the literature by using classical method for dealing with free boundary problems, cf. \cite{GunPro, Sol1, Sol2}
  for instance. However, the proofs are hard. By reducing (7.8) into the initial value problem (7.2) of a differential equation
  in the Banach manifold $\mathfrak{M}$, we can use some abstract result on invariant and quasi-invariant parabolic differential
  equations on Banach manifold under Lie group actions to give new proofs of these results which are relatively simpler. We
  leave this for future work.

  The above method can also be applied to more complex free boundary problems containing parabolic partial differential
  equations. Consider first the following free boundary problem modeling tumor growth:
\begin{equation}
\left\{
\begin{array}{ll}
  c\partial_t\sigma=\Delta\sigma-f(\sigma),&\quad x\in\Omega(t),\,\,t>0,\\
  -\Delta p=g(\sigma),&\quad  x\in\Omega(t),\,\,t>0,\\
  -\partial_{\mathbf{n}}\sigma=h(\sigma),&\quad  x\in\partial\Omega(t),\,\,t>0,\\
  p=\gamma\kappa,&\quad  x\in\partial\Omega(t),\,\,t>0,\\
  V_{\mathbf{n}}=-\partial_{\mathbf{n}}p,&\quad  x\in\partial\Omega(t),\,\,t>0,\\
  \sigma(x,0)=\sigma_0(x),&\quad  x\in\Omega_0,\\
  \Omega(0)=\Omega_{0},&
\end{array}\right.
\end{equation}
  where $c$ is a positive constant measuring the ratio between tumor cell doubling time scale and the nutrient diffusion
  time scale, and all other variables and constants are the same with those in (7.6). We assume that the function $f,g,h$
  satisfy the same conditions as in the problem (7.6).

  Let $\mathfrak{N}=\dot{\mathfrak{E}}^{m+\mu,m+\mu}({\mathbb{R}}^n)$ and
$$
  \mathfrak{N}_0=\{(\Omega,\sigma):\; (\Omega,\sigma)\in\dot{\mathfrak{E}}^{m+3+\mu,m+2+\mu}({\mathbb{R}}^n),\,
  \partial_{\mathbf{n}}\sigma+h(\sigma)=0\;\mbox{on}\;\partial\Omega\}.
$$
  $\mathfrak{N}_0$ is an embedded Banach submanifold of $\mathfrak{N}$. Given $(\Omega,\sigma)\in\mathfrak{N}_0$, let
  $p\in\dot{C}^{m+1+\mu}(\overline{\Omega})$ be the unique solution of the following linear elliptic boundary value problem:
$$
\left\{
\begin{array}{ll}
  -\Delta p=g(\sigma),&\quad \mbox{in}\;\;\Omega,\\
  p=\gamma\kappa_{\partial\Omega},&\quad \mbox{on}\;\;\partial\Omega.
\end{array}
\right.
$$
  We define a vector field $\mathcal{G}$ on $\mathfrak{N}$ with domain $\mathfrak{N}_0$ as follows:
$$
  \mathcal{G}(\Omega,\sigma)=(c^{-1}\Delta\sigma-c^{-1}f(\sigma),-\partial_{\mathbf{n}}p|_{\partial\Omega}), \quad
  \forall(\Omega,\sigma)\in\mathfrak{N}_0.
$$
  It follows that the problem (7.9) reduces into the following initial value problem for a differential equation on the
  Banach manifold $\mathfrak{N}$:
$$
\left\{
\begin{array}{ll}
  (\Omega(t),\sigma(t))'=\mathcal{G}(\Omega(t),\sigma(t)),&\quad t>0,\\
  (\Omega(0),\sigma(0))=(\Omega_0,\sigma_0).&
\end{array}
\right.
$$
  This problem, however, is more difficult to study than the problem (7.2) due to nonlinearity of the boundary value
  condition for $\sigma$ as it leads to the more complex Banach manifold $\mathfrak{N}_0$ of the domain of $\mathcal{G}$.

  Finally let us consider the following free boundary problem modeling the motion of liquid drops:
\begin{equation}
\left\{
\begin{array}{cl}
   \partial_t\bfv=\Delta\bfv-(\bfv\cdot\nabla)\bfv-\nabla\pi, &\quad\;\; x\in\Omega(t), \;\; t>0,\\
   \nabla\cdot\bfv=0,   &\quad \;\; x\in\Omega(t), \;\; t>0,\\
   T(\bfv,\pi)\bfn=-\gamma\kappa\bfn,   &\quad \;\; x\in\partial\Omega(t), \;\; t>0,\\
   V_n=\bfv\cdot\bfn, &\quad \;\; x\in\partial\Omega(t), \;\; t>0,\\
   \bfv|_{t=0}=\bfv_0, &\quad \;\; x\in\Omega_0,\\
   \Omega(0)=\Omega_0. &
\end{array}
\right.
\end{equation}
  Let
$$
  \mathfrak{N}==\{(\Omega,\bfv):\; (\Omega,\bfv)\in\vec{\dot{\mathfrak{E}}}^{m+\mu,m+\mu}({\mathbb{R}}^n),\,
  \nabla\cdot\bfv=0\;\mbox{in}\;\Omega\},
$$
$$
\begin{array}{rl}
  \mathfrak{N}_0=\{(\Omega,\bfv):\; (\Omega,\bfv)\in&\vec{\dot{\mathfrak{E}}}^{m+1+\mu,m+1+\mu}({\mathbb{R}}^n),\,
  \nabla\cdot\bfv=0\;\mbox{in}\;\Omega\\
  &\mbox{and}\;T(\bfv,\pi)\bfn+\gamma\kappa\bfn=0\;\mbox{on}\;\partial\Omega\}.
\end{array}
$$
  Here $\vec{\dot{\mathfrak{E}}}^{m+\mu,m+\mu}({\mathbb{R}}^n)$ denotes the Banach vector bundle over
  $\dot{\mathfrak{D}}^{m+\mu}({\mathbb{R}}^n)$ with fibre $[\dot{C}^{m+\mu}(\overline{\Omega})]^3$ for $\Omega\in
  \dot{\mathfrak{D}}^{m+\mu}({\mathbb{R}}^n)$, and similarly for $\vec{\dot{\mathfrak{E}}}^{m+1+\mu,m+1+\mu}({\mathbb{R}}^n)$.
  $\mathfrak{N}_0$ is an embedded Banach submanifold of $\mathfrak{N}$.
  We define a vector field $\mathcal{G}$ on $\mathfrak{N}$ with domain $\mathfrak{N}_0$ as follows:
$$
  \mathcal{G}(\Omega,\bfv)=(\Delta\bfv-\mathbb{P}[(\bfv\cdot\nabla)\bfv],\bfv\cdot\bfn), \quad
  \forall(\Omega,\bfv)\in\mathfrak{N}_0.
$$
  where $\mathbb{P}=I+\nabla(-\Delta)^{-1}\nabla$ is the denotes the Helmholtz-Weyl projection operator.
  Then the problem (7.10) reduces into an initial value problem for a differential equation on the
  Banach manifold $\mathfrak{N}$ of the following form:
$$
\left\{
\begin{array}{ll}
  (\Omega(t),\bfv(t))'=\mathcal{G}(\Omega(t),\bfv(t)),&\quad t>0,\\
  (\Omega(0),\bfv(0))=(\Omega_0,\bfv_0).&
\end{array}
\right.
$$

  In summary, by introducing the concepts of {\em differentiable points} in the Banach manifold $\dot{\mathfrak{D}}^{m+\mu}({\mathbb{R}}^n)$
  and {\em tangent spaces} of $\dot{\mathfrak{D}}^{m+\mu}({\mathbb{R}}^n)$ at such points, we see that evolutionary type free boundary
  problems can be nicely represented as differential equations in either $\dot{\mathfrak{D}}^{m+\mu}({\mathbb{R}}^n)$ or certain
  vector bundles over $\dot{\mathfrak{D}}^{m+\mu}({\mathbb{R}}^n)$, so that it becomes possible to treat such free boundary problems
  globally. Surely, in order to carry out this idea successfully, the structure and analytic properties of $\dot{\mathfrak{D}}^{m+\mu}({\mathbb{R}}^n)$
  as well as a well-developed theory of differential calculus in this Banach manifold should be systematically established. The author
  of the present paper is glad to invite interested reader to contribute intelligence and wisdom on this topic and its extensions
  such as the Banach manifold of simple $W^{m,p}$-domains in $\mathbb{R}^n$ ($m\in\mathbb{N}$, $1\leqslant p<\infty$),
  bounded $C^{m+\mu}$-domains in $\mathbb{R}^n$ of genus $N$ ($m\in\mathbb{N}$, $1\leqslant\mu\leqslant 1$,
  $N=1,2,\cdots$) and  bounded $W^{m,p}$-domains in $\mathbb{R}^n$ of genus $N$ ($m\in\mathbb{N}$, $1\leqslant p<\infty$,
  $N=1,2,\cdots$), etc. Such study clearly has extensive and strongly meaningful potential applications.


\end{document}